\documentclass[reqno,11pt]{amsart}
\usepackage{amsmath,amssymb,latexsym,cancel,rotating,hyperref}

\usepackage{eucal}

\usepackage{amscd}
\usepackage[dvips]{color}
\usepackage{multicol}
\usepackage[all]{xy}           
\usepackage{graphicx}
\usepackage{color}
\usepackage{colordvi}
\usepackage{xspace}
\usepackage{tikz}

\usepackage{enumitem}

\numberwithin{equation}{section}

\newtheorem{theorem}{Theorem}[section]
\newtheorem{proposition}[theorem]{Proposition}

\newtheorem{corollary}[theorem]{Corollary}
\newtheorem{lemma}[theorem]{Lemma}

\newtheorem{definition}[theorem]{Definition}

\input xy
\xyoption{poly}
\xyoption{2cell}
\xyoption{all}

\def\e{{\bf{e}}}

\def\b{{\bf{b}}}

\renewcommand{\eqref}[1]{{\rm (\ref{#1})}}

\begin{document}

\title[From the  quantum group $U_{q^{1/2}}( {\widehat {\mathfrak{sl}}}_2 )$ to the Kronecker quantum cluster algebra]
{From the  quantum group $U_{q^{1/2}}( {\widehat {\mathfrak{sl}}}_2)$ to the Kronecker quantum cluster algebra}

\author{Junyuan Huang, Xueqing Chen, Ming Ding and Fan Xu}
\address{School of Mathematics and Information Science\\
Guangzhou University, Guangzhou 510006, P.R.China}
\email{jy4545@e.gzhu.edu.cn (J.Huang)}
\address{Department of Mathematics,
 University of Wisconsin-Whitewater\\
800 W. Main Street, Whitewater, WI.53190. USA}
\email{chenx@uww.edu (X.Chen)}
\address{School of Mathematics and Information Science\\
Guangzhou University, Guangzhou 510006, P.R.China}
\email{dingming@gzhu.edu.cn (M.Ding)}
\address{Department of Mathematical Sciences\\
Tsinghua University\\
Beijing 100084, P.~R.~China} \email{fanxu@mail.tsinghua.edu.cn (F.Xu)}




\dedicatory{Dedicated to Professor Yanan Lin on the Occasion of His 70th Birthday}

\maketitle

\begin{abstract}
Let $\mathcal A_q$ be the Kronecker quantum cluster algebra with principal coefficients. We construct an explicit algebra homomorphism from the whole quantum group ${U_{{q^{1/2}}}}( {{{\widehat {\mathfrak{sl}}}_2}} )$ to $\mathcal A_q$, whose images together with either of the two frozen variables generate $\mathcal A_q$. As an application, we obtain  algebra homomorphisms from the split $\imath$quantum groups of affine $\widehat {\mathfrak{sl}}_2$ to $\mathcal A_q$. Moreover, we explicitly determine the images of the positive and negative real and imaginary root vectors of ${U_{{q^{1/2}}}}( {{{\widehat {\mathfrak{sl}}}_2}} )$ in terms of quantum cluster variables.
\end{abstract}


\section{Background}

Cluster algebras were introduced by Fomin and Zelevinsky \cite{ca1,ca2} to study total positivity \cite{G3} and Lusztig's dual canonical bases in coordinate rings and their $q$-deformations\cite{G1,G2}.  Quantum deformations of  cluster algebras were  introduced  by Berenstein and Zelevinsky \cite{AA}, in which each cluster generates a  quantum torus algebra and mutations are defined by quantum exchange relations.

 Quantum groups were introduced by Drinfeld \cite{Dri85} and Jimbo \cite{Jim85} in 1985, which are
deformations of the universal enveloping algebras. Their relations with quantum cluster algebras have been studied through algebra homomorphisms. Schrader and Shapiro \cite{SS19} constructed a cluster realization of quantum groups of finite type $A$. They obtain an embedding  of the quantum group $U_{q}(sl_{n+1})$
into certain quantum cluster algebra. Note that, the embedding is not surjective and
the image  is also not known. For each simple Lie algebra $\mathfrak{g}$, Ivan Ip \cite{Ip18}  constructed an algebra embedding of the
quantum group $U_{q}(\mathfrak{g})$ into a quantum torus algebra associated with a certain quiver. However, it is not clear whether the image  is in the quantum cluster
algebra associated with this quiver. In \cite{Shen}, Shen provided a rigid cluster model to realize the quantum group $U_{q}(\mathfrak{g})$ for $\mathfrak{g}$ of finite types $ADE$.
Recently, Huang, Chen, Ding and Xu \cite{HJY} obtained the fundamental relations in  arbitrary quantum cluster algebras with principal coefficients and  deduced an algebra homomorphism from the (untwisted)  positive parts of quantum groups to the corresponding quantum cluster algebras. This homomorphism was also studied by Berenstein and Rupel \cite{BR},  Fu and Zhang \cite{FZ25} from the perspective of algebraic representation theory.

Quantum symmetric pairs were introduced by Letzter \cite{Let99} and extended to the Kac-Moody setting by Kolb \cite{Ko14}. Their coideal subalgebras are called $\imath$quantum groups. Universal split $\imath$quantum groups and their presentations were also studied (see \cite{LW22,LW23,CLW21a}). Song \cite{So24} constructed cluster realizations of $\imath$quantum groups of finite type $\mathrm{AI}$. Lu, Ruan and Zhang \cite{LRZ25} constructed an analogue of Feigin's homomorphism from a universal split $\imath$quantum group to a quantum torus algebra. In affine type $A_1^{(1)}$, distinguished parameters give the $q^{1/2}$-Onsager algebra. Goff \cite{Goff24} constructed a homomorphism from the $q$-Onsager algebra to a quantum torus algebra and studied the images of its PBW elements in \cite{Goff25}.

In this paper, we focus on studying the relations between  the whole quantum group $U_{q^{1/2}}( {\widehat {\mathfrak{sl}}}_2 )$ and the Kronecker quantum cluster algebra  with principal coefficients. The paper is organized as follows. Section $2$ reviews the definitions of quantum cluster algebras and quantum groups, and fixes some conventions.
In Section $3$,  we establish an explicit algebra homomorphism $\Psi:U_{q^{1/2}}( {{{\widehat {\mathfrak{sl}}}_2}}  )\longrightarrow\mathcal A_q$, where $\mathcal A_q$ is the Kronecker quantum cluster algebra with principal coefficients. The images generate $A_q$ after adding either of the two frozen variables.  We also prove that the restriction of $\Psi$ to a rank-one quantum subalgebra is injective and identify its image  with $U_{q^{1/2}}( {{{ {\mathfrak{sl}}}_2}} )$.
In Section $4$, as an application, we obtain homomorphisms from split $\imath$quantum groups of affine type $A_1^{(1)}$ with arbitrary nonzero parameters, which particularly give a homomorphism from the $q^{1/2}$-Onsager algebra to $\mathcal A_q$.
In Section $5$, we explicitly determine the images of the positive and negative real and imaginary root vectors of the affine quantum group $U_{q^{1/2}}( {{{\widehat {\mathfrak{sl}}}_2}} )$ under $\Psi$ in terms of  the quantum cluster variables.

\section{Preliminaries}
For the convenience of the reader, we recall some  definitions and notations  concerning quantum cluster algebras \cite{AA} and quantum groups.

\subsection{Quantum cluster algebras}
For a positive integers $s$, we denote
by   $[1,s]$  the set $\{1,2,\ldots,s\}$.

\begin{definition}
A compatible pair is  $(\widetilde{B},\Lambda)$ such that
\begin{enumerate}
		\item $\widetilde{B}=(b_{ij})$  is an $m\times n$ $(m\ge n)$ integer matrix with its upper $n\times n$ submatrix  $B$  being skew-symmetrizable;
		\item $\Lambda=(\lambda_{ij})$ is an $m\times m$  skew-symmetric integer matrix;
        \item $\widetilde{B}^{T}\Lambda=\begin{bmatrix} D & {\bf{0}} \end{bmatrix}$, where $D=\operatorname{diag}(d_1, d_2, \dots, d_n)$ is a diagonal matrix
with $d_i$  $(i\in [1,n])$ positive integers.
	\end{enumerate}
 \end{definition}

Let $q$ be a  formal parameter and its square root $q^{\frac{1}{2}}$.
The  quantum torus $\mathcal T=\mathcal T(\Lambda)$ is the
$\mathbb Q(q^{\frac{1}{2}})$-algebra generated by
$\{X^{\e}:\e\in\mathbb Z^m\}$ with multiplication given by
\begin{equation}\label{(ef)}
	X^{\e}X^{\bf f}=q^{\frac{\Lambda(\e,\bf f)}{2}}X^{\e+\bf f},
	\qquad \e,\bf f\in\mathbb Z^m,
\end{equation}
where the matrix product $\e^{T}\Lambda\bf f$  is denoted by $\Lambda(\e,\bf f)$. The quantum torus $\mathcal T$ is an Ore domain and we denote by $\mathcal{F}$ its
fraction skew field.

Denote by $\e_i$  the $i$-th unit vector in $\mathbb{Z}^m$ for any $ i\in [1,m]$ and $X_i$ by $X^{\e_i}$, then
for an integer vector ${\bf{a}}=(a_1, a_2, \dots, a_m)\in \mathbb{Z}^m$,  we have $$X^{\bf{a}}=q^{\frac{\sum_{l<k}a_ka_l\lambda_{kl}}{2}} X_1^{a_1}X_2^{a_2}\cdots X_m^{a_m}.$$

The $\mathbb{Q}$-linear bar-involution $-: \mathcal{T} \longrightarrow \mathcal{T}$  is defined by
\[
\overline{f(q^{\frac{1}{2}})X^{\mathbf a}}
=f(q^{-\frac{1}{2}})X^{\mathbf a},
\qquad
f(q^{\frac{1}{2}})\in \mathbb Q(q^{\frac{1}{2}}),\quad \mathbf {\mathbf a} \in\mathbb{Z}^m.
\]

\begin{definition}\label{def of generalized seed}
The initial quantum seed is the triple $\big(\widetilde{\mathbf{x}},\Lambda,\widetilde{B}\big)$, where $(\widetilde{B},\Lambda)$ is a compatible pair, $\widetilde{\mathbf{x}}=\{X_{1},X_2,\ldots, X_{m}\}$ is the extended cluster, $\mathbf{x}=\{X_{1},X_2,\ldots,X_{n}\}$ is the cluster,
 $X_i$ $(i\in [1,n])$ are  quantum cluster variables and  $X_{i}$ $(i\in [n+1, m])$ are  frozen variables.
\end{definition}

For an integer $a$, we use the notation $[a]_{+} := max(a, 0).$ For an integer vector ${\bf{a}}=(a_1, a_2, \dots, a_m)\in \mathbb{Z}^m$,  we denote by
${[\bf{a}]_{+}}=([a_1]_{+}, [a_2]_{+}, \dots, [a_m]_{+}).$

\begin{definition}
	The mutation of a quantum seed $\big(\widetilde{\mathbf{x}},\Lambda,\widetilde{B}\big)$ in  direction $k \in [1,n]$ is a new quantum seed \smash{$\mu_k\big(\widetilde{\mathbf{x}},\Lambda,\widetilde{B}\big):=\big(\widetilde{\mathbf{x}}',
		\Lambda',\widetilde{B}'\big)$} defined by the following rule:
	\begin{enumerate}\itemsep=0pt
		\item the new extended cluster $\widetilde{\mathbf{x}}':=(\widetilde{\mathbf{x}}- \{X_k\})\cup\{X'_{k}\}$ with
		\begin{equation}\label{(er)}
			X'_k=X^{-\e_k+[\b_k]_+}+X^{-\e_k+[-\b_k]_+},
		\end{equation}
      where $\b_k$ is the $k$-column of $\widetilde{B}$;
		\item
		the matrix $\widetilde{B}':=\mu_k\big(\widetilde{B}\big)$ with
		\begin{gather*}\label{(B)}
			b'_{ij}=
			\begin{cases}
				-b_{ij}, &\text{ if }  i=k \text{ or } j=k;
				\\
				b_{ij}+\displaystyle\frac{|b_{ik}|b_{kj}+b_{ik}|b_{kj}|}{2}, &\text{otherwise};
			\end{cases}
		\end{gather*}
		\item
		the skew-symmetric matrix $\Lambda':=\mu_k\big(\Lambda\big)$ with
		\begin{gather*}\label{(A)}
			\lambda'_{ij}=
			\begin{cases}
				\lambda_{ij}, &\text{ if } i,j\ne k;
				\\
				- {\lambda _{ij}} + \sum\limits_{t = 1}^m {{{[{b_{ti}}]}_ + }{\lambda _{tj}}} , &\text{ if } i= k, j\ne k.
			\end{cases}
		\end{gather*}
	\end{enumerate}
\end{definition}

It is known that $\mu_k$ is an involution.  Two quantum seeds are called mutation-equivalent if they can be obtained from
each other through a sequence of mutations. Denote
$$\mathop{\mathbb{QP}}:=\mathbb Q(q^{\frac{1}{2}})[X^{\pm}_{n+1},\ldots,X^{\pm}_m].$$

\begin{definition}\label{def of gca}
	The quantum cluster algebra  $\mathcal A_{q} =\smash{\mathcal{A}\big(\widetilde{\mathbf{x}}, \Lambda, \widetilde{B}\big)}$ associated to the initial quantum seed $\big(\widetilde{\mathbf{x}}, \Lambda, \widetilde{B}\big)$ is the $\mathop{\mathbb{QP}}$-subalgebra of $\mathcal{F}$ generated by all quantum cluster variables from all quantum seeds mutation-equivalent to $\smash{\big(\widetilde{\mathbf{x}},\Lambda,\widetilde{B}\big)}$.
\end{definition}

\begin{definition}
	A  quantum cluster algebra $\smash{\mathcal{A}\big(\widetilde{\mathbf{x}}, \Lambda, \widetilde{B}\big)}$ is called principal coefficients quantum cluster algebra if the submatrix composed of the last $m-n$ rows of $\widetilde{B}$  is an identity matrix.
\end{definition}

\subsection{Quantum groups}
Let $\mathfrak g$ be a Kac-Moody Lie algebra  with rank $n$ and  $C = (c_{ij})_{1 \le i, j \le n}$  the corresponding Cartan
matrix with a symmetrizer $D=\operatorname{diag}(d_1,\ldots,d_n)$.  For any positive integer $m$, define

\begin{eqnarray*}
[m]_{q_i} = \frac{q_i^{m} - q_i^{-m}}{q_i - q_i^{-1}},\quad
{}[m]_{q_i}^! = \prod_{s=1}^{m} [s]_{q_i},\quad
{}\begin{bmatrix} m \\ r \end{bmatrix}_{q_i} = \frac{[m]_{q_i}^!}{[r]_{q_i}^! [m - r]_{q_i}^!},
\end{eqnarray*}
where we denote $q^{\frac{d_i}{2}}$ by $q_i$ for any $i\in [1,n]$ and $[0]_{q_i}^!=1$.

\begin{definition}
	 The quantum group $U_{q^{1/2}}(\mathfrak g)$ associated with $C$ is the $\mathbb Q(q^{\frac{1}{2}})$-algebra with generators
$E_i, F_i, K_i^{\pm1}$  $(1 \le i \le n)$	and the following  relations:
\begin{eqnarray}
&&K_i K_i^{-1} = K_i^{-1} K_i = 1; \hspace{0.2cm} K_i K_j = K_j K_i;\label{(cartan1)}\\
&&K_i E_j = q_i^{c_{ij}} E_j K_i;\hspace{0.2cm} K_i F_j = q_i^{-c_{ij}} F_j K_i;\label{(cartan2)}\\
&&[E_i, F_i] = \frac{K_i - K_i^{-1}}{q_i - q_i^{-1}};\hspace{0.2cm} [E_i, F_j] = 0,\hspace{0.2cm} i \ne j;\label{(cartan3)}
\end{eqnarray}	
and for  $i \ne j$,  the quantum Serre relations	
\begin{eqnarray}
	\sum_{r=0}^{1-c_{i j}}(-1)^r\left[\begin{array}{c}1-c_{i j} \\ r\end{array}\right]_{q_i} E_i^{1-c_{i j}-r} E_j E_i^r=0;\label{(2.6)}\\
	\sum_{r=0}^{1-c_{i j}}(-1)^r\left[\begin{array}{c}1-c_{i j} \\ r\end{array}\right]_{q_i} F_i^{1-c_{i j}-r} F_j F_i^r=0,\label{(2.7)}
\end{eqnarray}
where $[a,b]:=ab-ba.$
\end{definition}

 \section{An algebra homomorphism from  $U_{q^{1/2}}( {\widehat {\mathfrak{sl}}}_2 )$ to $\mathcal A_q$}
In this section, we mainly construct an algebra homomorphism from the  quantum group of affine type $A_1^{(1)}$ to the Kronecker quantum cluster algebra with principal coefficients, and then study its image.

We consider the following compatibility pair $(\widetilde B,\Lambda)$, where

$$\widetilde{B}=\binom{B}{I_2}=\left(\begin{array}{cc}0 & 2 \\ -2 & 0 \\ 1 & 0 \\ 0 & 1\end{array}\right);  \Lambda=\left(\begin{array}{cccc}0 & 0 & -1 & 0 \\ 0 & 0 & 0 & -1 \\ 1 & 0 & 0 & -2 \\ 0 & 1 & 2 & 0\end{array}\right).$$
It is easy to see that
$$\widetilde{B}^T \Lambda=(D \mid 0),$$
where $D$ is the $2\times 2$ identity matrix.
The corresponding Cartan matrix is
$$C=\left(\begin{array}{cc}2 & -2 \\ -2 & 2\end{array}\right).$$
Note that, $d_1 = d_2 = 1$ and $q_1 = q_2 = q^{\frac{1}{2}}$.

We set
$$x_1=X^{(1,0,0,0)};\hspace{0.2cm} x_2=X^{(0,1,0,0)};\hspace{0.2cm} y_1=X^{(0,0,1,0)};\hspace{0.2cm} y_2=X^{(0,0,0,1)}.$$

Let $\widetilde{\mathbf{x}}=\left(x_1, x_2, y_1, y_2\right)$, then $\Sigma=(\widetilde{\mathbf{x}}, \Lambda, \widetilde{B})$ is an initial quantum seed. Note that $y_1, y_2$ are frozen variables. The quantum cluster algebra  $\mathcal A_q$ associated with the initial quantum seed $\Sigma$ is called the Kronecker quantum
cluster algebra with principal coefficients.

In \cite{DXC23}, all cluster are denoted by $\{x_m,x_{m+1}| m \in \mathbb{Z}\}$ and thus all cluster variables are  $\{x_m| m \in \mathbb{Z}\}$.
In $\mathcal A_q$, there is an important  element defined by
$$X_\delta=X^{(-1,-1,1,0)}+X^{(-1,1,0,0)}+X^{(1,-1,1,1)}.$$

The $n$-th Chebyshev polynomials of the second kind associated to $X_\delta$ are defined recursively by $S_0=1, S_1=X_\delta$ and $S_{n+1}=S_n X_\delta-X^{(0,0,1,1)} S_{n-1}$ for any $n \ge 1$. Note that $\overline{S_n}=S_n$ for any $n \ge 0$.

\begin{theorem}\cite{DXC23}\label{3.2}
	In $\mathcal{A}_q$,  we have
	\begin{flalign*}
		&\quad(1) \hspace{0.2cm} X_\delta x_0  =x_{-1}+q^{-\frac{1}{2}} x_1 y_2; &&\\
		&\quad(2) \hspace{0.2cm} X_\delta x_{-m}  =x_{-(m+1)}+q y_1 y_2 x_{-(m-1)}, \hspace{0.2cm}\text{for $m \ge 1$;} &&\\
		&\quad(3) \hspace{0.2cm} x_1 X_\delta  =x_2+q^{-\frac{1}{2}} x_0 y_1; &&\\
		&\quad(4) \hspace{0.2cm} x_m X_\delta  =x_{m+1}+q x_{m-1} y_1 y_2,\hspace{0.2cm}\text{for $m \ge 2$} .&&
	\end{flalign*}
\end{theorem}
\begin{theorem}\cite{DXC23}\label{3.1}
	In  $\mathcal{A}_q$, we have
\begin{flalign*}
	&\quad(1) \hspace{0.2cm} x_{-(m+n+2)} x_{-m}=q^{-\frac{1}{2}} S_n X^{(0,0, m, m+1)}+x_{-(m+n+1)} x_{-(m+1)},\hspace{0.2cm}\text{for $m \ge 0$ and $n \ge 0$;} &&\\
	&\quad(2) \hspace{0.2cm} x_{m+n+2} x_m=q^{\frac{1}{2}} S_n X^{(0,0, m, m-1)}+x_{m+n+1} x_{m+1}, \hspace{0.2cm}\text{for $m \ge 1$ and $n \ge 0$;} &&\\
	&\quad(3) \hspace{0.2cm} x_0 x_{n+2}=q^{\frac{2 n+1}{2}} x_1 x_{n+1} y_2+S_n, \hspace{0.2cm}\text{for $n \ge 0$;} &&\\
	&\quad(4) \hspace{0.2cm} x_{-m} x_1=S_{m-1}+q^{\frac{1}{2}} x_{-(m-1)} x_0 y_1, \hspace{0.2cm}\text{for $m \ge 1$;} &&\\
	&\quad(5) \hspace{0.2cm} x_{-m} x_n=q x_{-(m-1)} x_{n-1} X^{(0,0,1,1)}+S_{m+n-2}, \hspace{0.2cm}\text{for $m \ge 1$ and $n \ge 2$.}&&
\end{flalign*}	
\end{theorem}

It is easy to show that
\begin{equation}\label{(3.8)}
	X^{(0,0, a, b)} X_\delta=q^{-a+b} X_\delta X^{(0,0, a, b)}.
\end{equation}

\begin{corollary}\label{c3.3}
	For any $m\ge 0$ and $a,b\in\mathbb Z$, we have
\begin{flalign}
	&\quad(1) \hspace{0.2cm}	x_{-m}y_1^a y_2^b=q^{am+b(1-m)}y_1^a y_2^b x_{-m};&&\label{(3.1)}\\
	&\quad(2) \hspace{0.2cm}	x_{m+3}y_1^a y_2^b=q^{a(m+1)-b(m+2)}y_1^a y_2^b x_{m+3}.&&\label{(3.2)}
\end{flalign}
\end{corollary}
\begin{proof}
	For any  $a,b\in\mathbb Z$, we prove both identities by induction on $m$. Note that, we have
	\begin{eqnarray*}
	&&x_1y_1^a y_2^b=q^{-a}y_1^a y_2^b x_1; \hspace{0.2cm}
		x_2y_1^a y_2^b=q^{-b}y_1^a y_2^b x_2;\\
		&&y_2y_1^a y_2^b=q^{2a}y_1^a y_2^b y_2;\hspace{0.2cm}
		y_1y_2y_1^a y_2^b=q^{2a-2b}y_1^a y_2^by_1y_2.\\
	\end{eqnarray*}
	
	When $m=0$, by (\ref{(er)}), we have
	\begin{eqnarray*}
		x_0=X^{(0,-1,0,0)}+X^{(2,-1,0,1)};\hspace{0.2cm} x_3=X^{(-1,0,1,0)}+X^{(-1,2,0,0)}.
	\end{eqnarray*}
Thus, we obtain
	\begin{eqnarray*}
		x_0y_1^a y_2^b=q^b y_1^a y_2^b x_0;\hspace{0.2cm} x_3y_1^a y_2^b=q^{a-2b}y_1^a y_2^b x_3.
	\end{eqnarray*}
	
When $m=1$, by Theorem~\ref{3.2}, we have
	$x_{-1}=X_\delta x_0-q^{-\frac12}x_1y_2$ and 	$x_4=x_3X_\delta-qx_2y_1y_2$. Hence, by (\ref{(3.8)}), we obtain
	\begin{eqnarray*}
		x_{-1}y_1^a y_2^b&=&\left(X_\delta x_0-q^{-\frac12}x_1y_2\right)y_1^a y_2^b=q^{a-b+b}y_1^a y_2^bX_\delta x_0-q^{-\frac12+2a-a}y_1^a y_2^b x_1y_2\\
		&=&q^a y_1^a y_2^b\left(X_\delta x_0-q^{-\frac12}x_1y_2\right)=q^a y_1^a y_2^b x_{-1};\\
		x_4y_1^a y_2^b&=&\left(x_3X_\delta-qx_2y_1y_2\right)y_1^a y_2^b=q^{a-b+a-2b}y_1^a y_2^b x_3X_\delta-q^{1+2a-2b-b}y_1^a y_2^b x_2y_1y_2\\
		&=&q^{2a-3b}y_1^a y_2^b\left(x_3X_\delta-qx_2y_1y_2\right)=q^{2a-3b}y_1^a y_2^b x_4.
	\end{eqnarray*}
	
	Suppose that this case holds for  $0\le m\le k+1$. When $m=k+2$,	by Theorem~\ref{3.2}, we have
	\begin{eqnarray*}
		x_{-(k+2)}y_1^a y_2^b&=&\left(X_\delta x_{-(k+1)}-qy_1y_2x_{-k}\right)y_1^a y_2^b\\
		&=&q^{a-b+a(k+1)-bk}y_1^a y_2^bX_\delta x_{-(k+1)}-q^{1+2a-2b+ak+b(1-k)}y_1^a y_2^b y_1y_2x_{-k}\\
		&=&q^{a(k+2)-b(k+1)}y_1^a y_2^b\left(X_\delta x_{-(k+1)}-qy_1y_2x_{-k}\right)\\
		&=&q^{a(k+2)-b(k+1)}y_1^a y_2^b x_{-(k+2)};\\
		x_{k+5}y_1^a y_2^b&=&\left(x_{k+4}X_\delta-qx_{k+3}y_1y_2\right)y_1^a y_2^b\\
		&=&q^{a-b+a(k+2)-b(k+3)}y_1^a y_2^b x_{k+4}X_\delta-q^{1+2a-2b+a(k+1)-b(k+2)}y_1^a y_2^b x_{k+3}y_1y_2\\
		&=&q^{a(k+3)-b(k+4)}y_1^a y_2^b\left(x_{k+4}X_\delta-qx_{k+3}y_1y_2\right)\\
&=&q^{a(k+3)-b(k+4)}y_1^a y_2^b x_{k+5}.
	\end{eqnarray*}

	This completes the proof.
\end{proof}

To construct an algebra homomorphism
$\Psi: U_{q^{1/2}}(\widehat{\mathfrak{s l}}_2) \longrightarrow \mathcal{A}_q$, we first introduce the following elements

\begin{eqnarray}
	 w_1&:=&-\frac{q y_2^{-1} x_3+q^{-\frac{1}{2}} X^{(0,0,-1,1)}x_{-1}}{\left(q^{\frac{1}{2}}-q^{-\frac{1}{2}}\right)^2}; \label{(w1)}\\
	 w_2&:=&-\frac{q^{\frac{1}{2}} X^{(0,0,-1,-1)} x_4+q X^{(0,0,-1,2)}x_0}{\left(q^{\frac{1}{2}}-q^{-\frac{1}{2}}\right)^2} . \label{(w2)}
\end{eqnarray}
Note that $\overline{w_1}=w_1$ and $\overline{w_2}=w_2$.

We define a map on the generators of $U_{q^{\frac{1}{2}}}(\widehat{\mathfrak{s l}}_2)$  as follows:
\begin{eqnarray}
	 &&\Psi\left(E_1\right)=x_1;\hspace{0.2cm} \Psi\left(E_2\right)=w_2; \hspace{0.2cm}
	 \Psi\left(F_1\right)=w_1; \hspace{0.2cm} \Psi\left(F_2\right)=X^{(0,1,0,-1)}; \nonumber\\
	 &&\Psi\left(K_1\right)=\Psi\left(K_2^{-1}\right)=X^{(0,0,1,-1)}; \hspace{0.2cm}
	 \Psi\left(K_2\right)=\Psi\left(K_1^{-1}\right)=X^{(0,0,-1,1)} .\label{(3.6)}
\end{eqnarray}

Note that the  images of generators are built from initial, projective and injective cluster variables and frozen monomials.
\begin{lemma}\label{3.3}
 The images in (\ref{(3.6)}) satisfy all  relations (\ref{(cartan1)}),(\ref{(cartan2)}) and (\ref{(cartan3)}) of $U_{q^{1/2}}(\widehat{\mathfrak{s l}}_2)$.
\end{lemma}
\begin{proof}
 First, we have
$$
\Psi\left(K_1\right) \Psi\left(K_2\right)=X^{(0,0,1,-1)} X^{(0,0,-1,1)}=X^{(0,0,-1,1)} X^{(0,0,1,-1)}=\Psi\left(K_2\right) \Psi\left(K_1\right)=1 .
$$
For $\mathbf{a}=\left(a_1, a_2, a_3, a_4\right)$, we have $\Lambda\left(\mathbf{e}_3-\mathbf{e}_4, \mathbf{a}\right)=a_1-a_2-2 a_3-2 a_4$. Since the quasicommutation factor in (\ref{(ef)}) is $q^{\Lambda\left(\mathbf{e}_3-\mathbf{e}_4, \mathbf{a}\right)}$, substituting the exponents in (\ref{(w1)}) and (\ref{(w2)}) gives
\begin{eqnarray*}
	&&\Psi\left(K_1\right) \Psi\left(E_1\right)=q \Psi\left(E_1\right) \Psi\left(K_1\right); \hspace{0.2cm}
	 \Psi\left(K_1\right) \Psi\left(E_2\right)=q^{-1} \Psi\left(E_2\right) \Psi\left(K_1\right); \\
	&&\Psi\left(K_1\right) \Psi\left(F_1\right)=q^{-1} \Psi\left(F_1\right) \Psi\left(K_1\right); \hspace{0.2cm}
	 \Psi\left(K_1\right) \Psi\left(F_2\right)=q \Psi\left(F_2\right) \Psi\left(K_1\right) .
\end{eqnarray*}
Since $\Psi\left(K_2\right)=\Psi\left(K_1\right)^{-1}$, the remaining  relations of (\ref{(cartan1)}) and (\ref{(cartan2)}) follow.

Note that, by Theorem \ref{3.1} and applying the bar-involution, we have
\begin{eqnarray*}
	{\left[x_1, q y_2^{-1} x_3\right] } & =&-\left(q^{\frac{1}{2}}-q^{-\frac{1}{2}}\right) X^{(0,0,1,-1)}; \\
	{\left[x_1, q^{-\frac{1}{2}} X^{(0,0,-1,1)} x_{-1}\right] } & =&\left(q^{\frac{1}{2}}-q^{-\frac{1}{2}}\right) X^{(0,0,-1,1)}; \\
	{\left[q^{\frac{1}{2}} X^{(0,0,-1,-1)} x_4, X^{(0,1,0,-1)}\right] } & =&\left(q^{\frac{1}{2}}-q^{-\frac{1}{2}}\right) X^{(0,0,1,-1)}; \\
	{\left[q X^{(0,0,-1,2)} x_0, X^{(0,1,0,-1)}\right] } & =&-\left(q^{\frac{1}{2}}-q^{-\frac{1}{2}}\right) X^{(0,0,-1,1)} .
\end{eqnarray*}
Substituting these identities into (\ref{(w1)}) and (\ref{(w2)}), we obtain
\begin{eqnarray*}
\left[\Psi\left(E_1\right), \Psi\left(F_1\right)\right]=\frac{X^{(0,0,1,-1)}-X^{(0,0,-1,1)}}{q^{\frac{1}{2}}-q^{-\frac{1}{2}}}=\frac{\Psi\left(K_1\right)-\Psi\left(K_1\right)^{-1}}{q_1-q_1^{-1}};\\
\left[\Psi\left(E_2\right), \Psi\left(F_2\right)\right]=\frac{X^{(0,0,-1,1)}-X^{(0,0,1,-1)}}{q^{\frac{1}{2}}-q^{-\frac{1}{2}}}=\frac{\Psi\left(K_2\right)-\Psi\left(K_2\right)^{-1}}{q_2-q_2^{-1}} .
\end{eqnarray*}
Since $\Lambda\left(\mathbf{e}_1, \mathbf{e}_2-\mathbf{e}_4\right)=0$, we have $\left[\Psi\left(E_1\right), \Psi\left(F_2\right)\right]=0$. To prove $\left[\Psi\left(E_2\right), \Psi\left(F_1\right)\right]=0$, first we have $x_3 x_4=x_4 x_3$ and $x_{-1} x_0=x_0 x_{-1}$. By (2.1), we obtain
\begin{eqnarray*}
	{\left[q y_2^{-1} x_3, q^{\frac{1}{2}} X^{(0,0,-1,-1)} x_4\right] } & =&0; \\
	{\left[q^{-\frac{1}{2}} X^{(0,0,-1,1)} x_{-1}, q X^{(0,0,-1,2)} x_0\right] } & =&0; \\
	{\left[q y_2^{-1} x_3, q X^{(0,0,-1,2)} x_0\right] } & =&X^{(0,0,-1,1)}\left(q^{-2} x_3 x_0-x_0 x_3\right); \\
	{\left[q^{-\frac{1}{2}} X^{(0,0,-1,1)} x_{-1}, q^{\frac{1}{2}} X^{(0,0,-1,-1)} x_4\right] } & =&q^{-3} y_1^{-2}\left(x_{-1} x_4-x_4 x_{-1}\right) .
\end{eqnarray*}
By Theorem \ref{3.2} and the bar-involution, we have $x_{-1} x_4-x_4 x_{-1}=X^{(0,0,1,1)}\left(q x_0 x_3-\right.$ $q^{-1} x_3 x_0$ ). Note that $q^{-3} y_1^{-2} X^{(0,0,1,1)}=q^{-1} X^{(0,0,-1,1)}$.

Therefore
\begin{eqnarray*}
	&&{\left[w_1, w_2\right]}\\&= & \frac{1}{\left(q^{\frac{1}{2}}-q^{-\frac{1}{2}}\right)^4}\left(\left[q y_2^{-1} x_3, q^{\frac{1}{2}} X^{(0,0,-1,-1)} x_4\right]+\left[q y_2^{-1} x_3, q X^{(0,0,-1,2)} x_0\right]\right. \\
&	& \left.+\left[q^{-\frac{1}{2}} X^{(0,0,-1,1)} x_{-1}, q^{\frac{1}{2}} X^{(0,0,-1,-1)} x_4\right]+\left[q^{-\frac{1}{2}} X^{(0,0,-1,1)} x_{-1}, q X^{(0,0,-1,2)} x_0\right]\right) \\
&	= & 0 .
\end{eqnarray*}
Hence $\left[\Psi\left(E_2\right), \Psi\left(F_1\right)\right]=0$. The proof is complete.
\end{proof}

\begin{lemma}
 In  $\mathcal{A}_q$, we have
\begin{flalign}
	&\quad(1)\hspace{0.2cm} w_2 x_1-q^{-1} x_1 w_2  =-\frac{q\left(1-q^{-2}\right)}{\left(q^{\frac{1}{2}}-q^{-\frac{1}{2}}\right)^2} y_2^{-1} X_\delta;  &&\label{(w2x1-x1w2)}\\
	&\quad(2)\hspace{0.2cm} w_1 X^{(0,1,0,-1)}-q X^{(0,1,0,-1)} w_1  =\frac{q^{-1}\left(q^2-1\right)}{\left(q^{\frac{1}{2}}-q^{-\frac{1}{2}}\right)^2} y_1^{-1} X_\delta;  &&\label{(w1F2-F2w1)}\\
	&\quad(3)\hspace{0.2cm} w_1 y_1^{-1} X_\delta-q^{-1} y_1^{-1} X_\delta w_1  =\left(q^{\frac{1}{2}}-q^{-\frac{1}{2}}\right) w_2; &&\label{(w1xd-xdw1)}\\
	&\quad(4)\hspace{0.2cm} {\left[w_2, q^{-\frac{1}{2}} y_1^{-1} X_\delta\right] }  =-\left(1-q^{-1}\right) X^{(0,0,-1,1)} w_1;  &&\label{(w2xd-xdw2)}\\
	&\quad(5)\hspace{0.2cm} x_1 X_\delta-q^{-1} X_\delta x_1  =\left(1-q^{-1}\right) x_2;  &&\label{(x1xd-xdx1)}\\
	&\quad(6)\hspace{0.2cm} X^{(0,1,0,-1)} y_1^{-1} X_\delta-q y_1^{-1} X_\delta X^{(0,1,0,-1)}  =q^{\frac{1}{2}}(1-q) x_1 . &&\label{(F2xd-xdF2)}
\end{flalign}
\end{lemma}
\begin{proof}
 By Theorem \ref{3.1} and the bar-involution, we obtain
\begin{eqnarray*}
	x_4 x_1-x_1 x_4  =\left(q^{\frac{3}{2}}-q^{-\frac{1}{2}}\right) y_1 X_\delta; \hspace{0.2cm}
	x_{-1} x_2-q^2 x_2 x_{-1}  =\left(1-q^2\right) X_\delta .
\end{eqnarray*}
 Substituting (\ref{(w1)}) and (\ref{(w2)}) and by (\ref{(ef)}), we have
\begin{eqnarray*}
	w_2 x_1-q^{-1} x_1 w_2 & =&-\frac{q^{\frac{1}{2}} X^{(0,0,-1,-1)}\left(x_4 x_1-x_1 x_4\right)}{\left(q^{\frac{1}{2}}-q^{-\frac{1}{2}}\right)^2}=-\frac{q-q^{-1}}{\left(q^{\frac{1}{2}}-q^{-\frac{1}{2}}\right)^2} y_2^{-1} X_\delta; \\
	w_1 X^{(0,1,0,-1)}-q X^{(0,1,0,-1)} w_1 & =&-\frac{q^{-1} y_1^{-1}\left(x_{-1} x_2-q^2 x_2 x_{-1}\right)}{\left(q^{\frac{1}{2}}-q^{-\frac{1}{2}}\right)^2}=\frac{q-q^{-1}}{\left(q^{\frac{1}{2}}-q^{-\frac{1}{2}}\right)^2} y_1^{-1} X_\delta .
\end{eqnarray*}
This proves (\ref{(w2x1-x1w2)}) and (\ref{(w1F2-F2w1)}).

By Theorem \ref{3.2} and the bar-involution, we obtain
\begin{eqnarray*}
&&x_3 X_\delta-q^{-1} X_\delta x_3  =\left(1-q^{-1}\right) x_4;\hspace{0.2cm}
	x_{-1} X_\delta-X_\delta x_{-1}  =(q-1) x_0 y_1 y_2; \\
&&	x_4 X_\delta-X_\delta x_4  =(1-q) y_1 y_2 x_3;\hspace{0.2cm}
	x_0 X_\delta-q X_\delta x_0  =(1-q) x_{-1} .
\end{eqnarray*}

Substituting these four identities into (\ref{(w1)}) and (\ref{(w2)}), and by (\ref{(3.8)}), we obtain
\begin{eqnarray*}
	w_1 y_1^{-1} X_\delta-q^{-1} y_1^{-1} X_\delta w_1 & =&-\frac{q-1}{\left(q^{\frac{1}{2}}-q^{-\frac{1}{2}}\right)^2}\left(X^{(0,0,-1,-1)} x_4+q^{\frac{1}{2}} X^{(0,0,-1,2)} x_0\right) \\
	& =&\left(q^{\frac{1}{2}}-q^{-\frac{1}{2}}\right) w_2;\\
	{\left[w_2, q^{-\frac{1}{2}} y_1^{-1} X_\delta\right] } & =&\frac{1-q^{-1}}{\left(q^{\frac{1}{2}}-q^{-\frac{1}{2}}\right)^2}\left(y_1^{-1} x_3+q^{-\frac{1}{2}} X^{(0,0,-2,2)} x_{-1}\right) \\
	& =&-\left(1-q^{-1}\right) X^{(0,0,-1,1)} w_1 .
\end{eqnarray*}
This proves (\ref{(w1xd-xdw1)}) and (\ref{(w2xd-xdw2)}).	
	
	Finally, by Theorem \ref{3.2} and the bar-involution, we have
		\begin{eqnarray*}
		x_1 X_\delta-q^{-1} X_\delta x_1  =\left(1-q^{-1}\right) x_2;\hspace{0.2cm}
		x_2 X_\delta-X_\delta x_2  =(1-q) y_1 y_2 x_1 .
	\end{eqnarray*}
	The first identity is (\ref{(x1xd-xdx1)}). Since $X^{(0,1,0,-1)}=q^{-\frac{1}{2}} x_2 y_2^{-1}$ and $X^{(0,0,-1,-1)} y_1 y_2=q^{-1}$, the second identity gives
	$$X^{(0,1,0,-1)} y_1^{-1} X_\delta-q y_1^{-1} X_\delta X^{(0,1,0,-1)}=q^{\frac{3}{2}} X^{(0,0,-1,-1)}\left(x_2 X_\delta-X_\delta x_2\right)=q^{\frac{1}{2}}(1-q) x_1 .$$
	Thus (\ref{(F2xd-xdF2)}) follows. This proof is complete.
\end{proof}

\begin{lemma}\label{3.5}
 The images in (\ref{(3.6)}) satisfy quantum Serre relations
\begin{flalign*}
&\quad(1)\hspace{0.2cm}	\sum_{r=0}^3(-1)^r\left[\begin{array}{l}
		3 \\
		r
	\end{array}\right]_{q^{\frac{1}{2}}} \Psi\left(E_i\right)^{3-r} \Psi\left(E_j\right) \Psi\left(E_i\right)^r=0 ;&&\\
	&\quad(2)\hspace{0.2cm} \sum_{r=0}^3(-1)^r\left[\begin{array}{l}
		3 \\
		r
	\end{array}\right]_{q^{\frac{1}{2}}} \Psi\left(F_i\right)^{3-r} \Psi\left(F_j\right) \Psi\left(F_i\right)^r=0.&&
\end{flalign*}
\end{lemma}
\begin{proof}
Note that $[3]_{q^{\frac{1}{2}}}=q+1+q^{-1}$.

 (1) We prove that the first equation for $i=2, j=1$.  We compute
 	\begingroup
 \interdisplaylinepenalty=0
\begin{eqnarray*}
	&& w_2^3 x_1-[3]_{q^{\frac{1}{2}}} w_2^2 x_1 w_2+[3]_{q^{\frac{1}{2}}} w_2 x_1 w_2^2-x_1 w_2^3 \\
	&= & w_2^2\left(w_2 x_1-q^{-1} x_1 w_2\right)-(1+q) w_2\left(w_2 x_1-q^{-1} x_1 w_2\right) w_2+q\left(w_2 x_1-q^{-1} x_1 w_2\right) w_2^2\\
	&\mathop  = \limits^{(\ref{(w2x1-x1w2)})} &-\frac{q^{\frac{1}{2}}\left(1-q^{-2}\right)}{\left(q^{\frac{1}{2}}-q^{-\frac{1}{2}}\right)^2}\left(q^{\frac{1}{2}} w_2^2 y_2^{-1} X_\delta-q^{\frac{1}{2}}(1+q) w_2 y_2^{-1} X_\delta w_2+q^{\frac{3}{2}} y_2^{-1} X_\delta w_2^2\right)  \\
	& =&-\frac{q^{\frac{3}{2}}\left(1-q^{-2}\right)}{\left(q^{\frac{1}{2}}-q^{-\frac{1}{2}}\right)^2} X^{(0,0,1,-1)}\left(q w_2^2 q^{-\frac{1}{2}} y_1^{-1} X_\delta-(1+q) w_2 q^{-\frac{1}{2}} y_1^{-1} X_\delta w_2+q^{-\frac{1}{2}} y_1^{-1} X_\delta w_2^2\right) \\
	& =&-\frac{q^{\frac{3}{2}}\left(1-q^{-2}\right)}{\left(q^{\frac{1}{2}}-q^{-\frac{1}{2}}\right)^2} X^{(0,0,1,-1)}\left(q w_2\left[w_2, q^{-\frac{1}{2}} y_1^{-1} X_\delta\right]-\left[w_2, q^{-\frac{1}{2}} y_1^{-1} X_\delta\right] w_2\right)\\
	& \mathop  = \limits^{(\ref{(w2xd-xdw2)})}&\frac{q\left(1-q^{-2}\right)}{q^{\frac{1}{2}}-q^{-\frac{1}{2}}} X^{(0,0,1,-1)}\left(q w_2 X^{(0,0,-1,1)} w_1-X^{(0,0,-1,1)} w_1 w_2\right)  \\
	& =&\frac{q\left(1-q^{-2}\right)}{q^{\frac{1}{2}}-q^{-\frac{1}{2}}}\left(w_2 w_1-w_1 w_2\right)=\frac{q\left(1-q^{-2}\right)}{q^{\frac{1}{2}}-q^{-\frac{1}{2}}}\left[\Psi\left(E_2\right), \Psi\left(F_1\right)\right]=0.
\end{eqnarray*}	
Now, we prove that the first equation for $i=1, j=2$. We compute
\begin{eqnarray*}
	&& x_1^3 w_2-[3]_{q^{\frac{1}{2}}} x_1^2 w_2 x_1+[3]_{q^{\frac{1}{2}}} x_1 w_2 x_1^2-w_2 x_1^3 \\
	&	= & x_1^2\left(x_1 w_2-q w_2 x_1\right)-\left(1+q^{-1}\right) x_1\left(x_1 w_2-q w_2 x_1\right) x_1 +q^{-1}\left(x_1 w_2-q w_2 x_1\right) x_1^2\\
	&\mathop  = \limits^{(\ref{(w2x1-x1w2)})}& \frac{q^2\left(1-q^{-2}\right)}{\left(q^{\frac{1}{2}}-q^{-\frac{1}{2}}\right)^2}\left(x_1^2 y_2^{-1} X_\delta-\left(1+q^{-1}\right) x_1 y_2^{-1} X_\delta x_1+q^{-1} y_2^{-1} X_\delta x_1^2\right)\\
	&=&\frac{q^2\left(1-q^{-2}\right)}{\left(q^{\frac{1}{2}}-q^{-\frac{1}{2}}\right)^2} y_2^{-1}\left(x_1\left(x_1 X_\delta-q^{-1} X_\delta x_1\right)-\left(x_1 X_\delta-q^{-1} X_\delta x_1\right) x_1\right)\\
	&\mathop  = \limits^{(\ref{(x1xd-xdx1)})}& \frac{q^{\frac{3}{2}}\left(1-q^{-2}\right)}{q^{\frac{1}{2}}-q^{-\frac{1}{2}}} y_2^{-1}\left(x_1 x_2-x_2 x_1\right)=0.
\end{eqnarray*}

(2) We prove that the second equation for $i=1, j=2$. We compute
\begin{eqnarray*}
	&& w_1^3 X^{(0,1,0,-1)}-[3]_{q^{\frac{1}{2}}} w_1^2 X^{(0,1,0,-1)} w_1+[3]_{q^{\frac{1}{2}}} w_1 X^{(0,1,0,-1)} w_1^2-X^{(0,1,0,-1)} w_1^3 \\
	& =&w_1^2\left(w_1 X^{(0,1,0,-1)}-q X^{(0,1,0,-1)} w_1\right)-\left(1+q^{-1}\right) w_1\left(w_1 X^{(0,1,0,-1)}-q X^{(0,1,0,-1)} w_1\right) \\
	&& \cdot w_1+q^{-1}\left(w_1 X^{(0,1,0,-1)}-q X^{(0,1,0,-1)} w_1\right) w_1^2\\
	&\mathop  = \limits^{(\ref{(w1F2-F2w1)})} &\frac{q^{-1}\left(q^2-1\right)}{\left(q^{\frac{1}{2}}-q^{-\frac{1}{2}}\right)^2}\left(w_1^2 y_1^{-1} X_\delta-\left(1+q^{-1}\right) w_1 y_1^{-1} X_\delta w_1+q^{-1} y_1^{-1} X_\delta w_1^2\right)\\
	&\mathop  = \limits^{(\ref{(w1xd-xdw1)})} & \frac{q^{-1}\left(q^2-1\right)}{\left(q^{\frac{1}{2}}-q^{-\frac{1}{2}}\right)^2}\left(w_1\left(w_1 y_1^{-1} X_\delta-q^{-1} y_1^{-1} X_\delta w_1\right)-\left(w_1 y_1^{-1} X_\delta-q^{-1} y_1^{-1} X_\delta w_1\right) w_1\right)\\
	&	=& \frac{q^{-1}\left(q^2-1\right)}{q^{\frac{1}{2}}-q^{-\frac{1}{2}}}\left(w_1 w_2-w_2 w_1\right)=\frac{q^{-1}\left(q^2-1\right)}{q^{\frac{1}{2}}-q^{-\frac{1}{2}}}\left[\Psi\left(F_1\right), \Psi\left(E_2\right)\right]=0.
\end{eqnarray*}
Now, we prove that the second equation for $i=2, j=1$. We compute
\begin{eqnarray*}
	&& \left(X^{(0,1,0,-1)}\right)^3 w_1-[3]_{q^{\frac{1}{2}}}\left(X^{(0,1,0,-1)}\right)^2 w_1 X^{(0,1,0,-1)} \\
	&& +[3]_{q^{\frac{1}{2}}} X^{(0,1,0,-1)} w_1\left(X^{(0,1,0,-1)}\right)^2-w_1\left(X^{(0,1,0,-1)}\right)^3 \\
	&= & \left(X^{(0,1,0,-1)}\right)^2\left(X^{(0,1,0,-1)} w_1-q^{-1} w_1 X^{(0,1,0,-1)}\right) \\
	&& -(1+q) X^{(0,1,0,-1)}\left(X^{(0,1,0,-1)} w_1-q^{-1} w_1 X^{(0,1,0,-1)}\right) X^{(0,1,0,-1)} \\
	&	& +q\left(X^{(0,1,0,-1)} w_1-q^{-1} w_1 X^{(0,1,0,-1)}\right)\left(X^{(0,1,0,-1)}\right)^2\\
	&\mathop  = \limits^{(\ref{(w1F2-F2w1)})}&-\frac{q^{-2}\left(q^2-1\right)}{\left(q^{\frac{1}{2}}-q^{-\frac{1}{2}}\right)^2}\left(\left(X^{(0,1,0,-1)}\right)^2 y_1^{-1} X_\delta\right.\\
	&& \left.-(1+q) X^{(0,1,0,-1)} y_1^{-1} X_\delta X^{(0,1,0,-1)}+q y_1^{-1} X_\delta\left(X^{(0,1,0,-1)}\right)^2\right) \\
	&\mathop  = \limits^{(\ref{(F2xd-xdF2)})} & -\frac{q^{-2}\left(q^2-1\right)}{\left(q^{\frac{1}{2}}-q^{-\frac{1}{2}}\right)^2}\left(X^{(0,1,0,-1)}\left(X^{(0,1,0,-1)} y_1^{-1} X_\delta-q y_1^{-1} X_\delta X^{(0,1,0,-1)}\right)\right. \\
	&	& \left.-\left(X^{(0,1,0,-1)} y_1^{-1} X_\delta-q y_1^{-1} X_\delta X^{(0,1,0,-1)}\right) X^{(0,1,0,-1)}\right)\\
	&	=&\frac{q^{-\frac{3}{2}}\left(1-q^2\right)(1-q)}{\left(q^{\frac{1}{2}}-q^{-\frac{1}{2}}\right)^2}\left(X^{(0,1,0,-1)} x_1-x_1 X^{(0,1,0,-1)}\right)\\
	&=&\frac{q^{-\frac{3}{2}}\left(1-q^2\right)(1-q)}{\left(q^{\frac{1}{2}}-q^{-\frac{1}{2}}\right)^2}\left[\Psi\left(F_2\right), \Psi\left(E_1\right)\right]=0 .
\end{eqnarray*}	\endgroup
The proof is complete.
\end{proof}

\begin{theorem}\label{3.6}
	 The map $\Psi: U_{q^{1/2}}(\widehat{\mathfrak{s l}}_2) \longrightarrow \mathcal{A}_q$ is an algebra homomorphism.
\end{theorem}
\begin{proof}
	 The proof follows from Lemma \ref{3.3} and Lemma \ref{3.5}.
\end{proof}

The next result says that the images together with either of the two frozen variables generate $\mathcal A_q$.

 \begin{theorem}\label{3.8}
 		We have $\mathbb{Q}(q^{\frac{1}{2}})\left\langle\operatorname{Im} \Psi, y_2\right\rangle=\mathbb{Q}(q^{\frac{1}{2}})\left\langle\operatorname{Im} \Psi, y_1\right\rangle=\mathcal{A}_q .$
\end{theorem}
\begin{proof}
	By (\ref{(3.6)}), we have $y_1=q \Psi\left(K_1\right) y_2, y_2=q y_1 \Psi\left(K_2\right)$. Thus, we obtain
$$\mathbb{Q}(q^{\frac{1}{2}})\left\langle\operatorname{Im} \Psi, y_2\right\rangle=\mathbb{Q}(q^{\frac{1}{2}})\left\langle\operatorname{Im} \Psi, y_1\right\rangle .
$$

We only need to prove $\mathbb{Q}(q^{\frac{1}{2}})\left\langle\operatorname{Im} \Psi, y_1\right\rangle=\mathcal{A}_q$. The $\subseteq$ inclusion is obvious. Thus, we need to prove the $\supseteq$ inclusion.  According to  \cite{AA}, we obtain that $x_0, x_1, x_2, x_3, y_1^{ \pm 1}, y_2^{ \pm 1}$ generate $\mathcal{A}_q$ since the Kronecker quiver is acyclic. It is enough to show that $ x_0, x_1, x_2, x_3, y_1^{ \pm 1}, y_2^{ \pm 1} \in \mathbb{Q}(q^{\frac{1}{2}})\left\langle\operatorname{Im} \Psi, y_1\right\rangle$.

By (\ref{(w1F2-F2w1)}), we obtain
$$X_\delta=\frac{q\left(q^{\frac{1}{2}}-q^{-\frac{1}{2}}\right)^2}{q^2-1} y_1 \Psi\left(F_1 F_2-q F_2 F_1\right) .
$$
Thus, $X_\delta \in \mathbb{Q}(q^{\frac{1}{2}})\left\langle\operatorname{Im} \Psi, y_1\right\rangle$.
Hence, we obtain
\begingroup
\interdisplaylinepenalty=0
\begin{eqnarray*}
	&&x_1 =\Psi\left(E_1\right);\hspace{0.2cm}
	x_2=q^{\frac{1}{2}} y_1 \Psi\left(K_2 F_2\right);\hspace{0.2cm} 	
    x_3 = {X_\delta }{x_2} - {q^2}y_1^2{X^{(0,0, - 1,1)}}{x_1}; \\
	&&y_1^{-1}  ={\frac{q-1}{q+1}} x_1 \Psi\left(F_1 F_2-q F_2 F_1\right) \Psi\left(K_2 F_2\right)-q^{\frac{1}{2}} x_1 X^{(0,0,-1,1)} x_1-q^{\frac{3}{2}}\left(\Psi\left(K_2 F_2\right)\right)^2; \\
	&&y_2^{-1}  =q^{-1} \Psi\left(K_1\right) y_1^{-1};\hspace{0.2cm}
	x_0  =q^{-\frac{1}{2}} y_1^{-1} X_\delta x_1-q^{-\frac{1}{2}} y_1^{-1} x_2.
	\end{eqnarray*}
\endgroup
Hence $ x_0, x_1, x_2, x_3,y_1^{ \pm 1}, y_2^{ \pm 1} \in \mathbb{Q}(q^{\frac{1}{2}})\left\langle\operatorname{Im} \Psi, y_1\right\rangle.$ The proof is complete.
\end{proof}

 Now we consider the subalgebra $U_1=\left\langle E_1, F_1, K_1, K_1^{-1}\right\rangle$ of $U_{q^{1/2}}(\widehat{\mathfrak{s l}}_2)$. Note that, there is a natural algebra isomorphism $\iota_1: U_{q^{1/2}}(\mathfrak{s} \mathfrak{l}_2) \longrightarrow U_1$ defined by $\iota_1(E)=E_1$, $\iota_1(F)=F_1, \iota_1(K)=K_1$ and $\iota_1\left(K^{-1}\right)=K_1^{-1}$. Set
$$\mathcal{B}_{q^{\frac{1}{2}}}:=\mathbb{Q}\left(q^{\frac{1}{2}}\right)\left\langle x_1, w_1, X^{(0,0,-1,1)}, X^{(0,0,1,-1)}\right\rangle \subseteq \mathcal{A}_q .$$

\begin{theorem}\label{3.9}
	The restriction $\rho=\left.\Psi\right|_{U_1}: U_1 \longrightarrow \mathcal{B}_{q^{\frac{1}{2}}}$ is an algebra isomorphism, where $\rho\left(E_1\right)=x_1; \rho\left(F_1\right)=w_1; \rho\left(K_1\right)=X^{(0,0,1,-1)}$ and $\rho\left(K_1^{-1}\right)=X^{(0,0,-1,1)}$. In particular, $\mathcal{B}_{q^{1/2}} \simeq U_{q^{1/2}}(\mathfrak{s l}_2)$.
\end{theorem}
\begin{proof}
 By Theorem \ref{3.6}, $\rho$ is an algebra homomorphism. By the definition of $\mathcal{B}_{q^{1/2}}$, it is surjective. It remains to prove that it is injective. By the PBW theorem, the elements $F_1^a K_1^b E_1^c\ (a, c \in \mathbb{N}, b \in \mathbb{Z})$ form a $\mathbb{Q}\left(q^{\frac{1}{2}}\right)$-basis of $U_1$.

Note that
\begin{eqnarray*}
	w_1&= &-\frac{1}{\left(q^{\frac{1}{2}}-q^{-\frac{1}{2}}\right)^2}\left(X^{(-1,0,1,-1)}+X^{(-1,2,0,-1)}+X^{(-1,0,-1,1)}\right. \\
	&& \left.+X^{(-1,-2,0,1)}+\left(q^{\frac{1}{2}}+q^{-\frac{1}{2}}\right) X^{(1,-2,0,2)}+X^{(3,-2,0,3)}\right) .
	\end{eqnarray*}
Thus, $X^{(-1,2,0,-1)}$ is the only monomial in $w_1$ with a positive power of $x_2$. Therefore the unique monomial of maximal $\mathbf{e}_2$-coordinate in $\rho\left(F_1^a K_1^b E_1^c\right)=w_1^a\left(X^{(0,0,1,-1)}\right)^b X^{(c, 0,0,0)}$ is, up to a nonzero scalar,
\begin{equation}
	X^{(-(a-c), 2 a, b,-(a+b))} . \label{(3.13)}
\end{equation}

For fixed $a$, distinct pairs $(b, c)$ give distinct monomials in (\ref{(3.13)}).
Suppose that $\sum_{a, b, c} $ $\lambda_{a, b, c} \rho\left(F_1^a K_1^b E_1^c\right)=0$ is a finite linear relation. Choose $a_0$ maximal such that some $\lambda_{a_0, b, c} \neq 0$. The terms with $a<a_0$ have $\mathbf{e}_2$-coordinate strictly less than $2 a_0$. The terms of $\mathbf{e}_2$-coordinate $2 a_0$ coming from $a=a_0$ are distinct for distinct pairs $(b, c)$.  Hence $\lambda_{a_0, b, c}=0$ for all $b, c$, which is a contradiction. Thus  $\rho$ is injective. The proof is complete.
\end{proof}

\section{An application to the split $\imath$quantum group of affine $\widehat {\mathfrak{sl}}_2$} Let $\boldsymbol{\varsigma}=\left(\varsigma_1, \varsigma_2\right)$ with $\varsigma_1, \varsigma_2 \in$ $\mathbb{Q}(q^{\frac{1}{2}}) \backslash\{0\}$.
\begin{definition}\label{4.1}
	 Following \cite{Let99,Ko14}, the split $\imath$quantum group $U_{\varsigma}^{\imath}$ is the $\mathbb{Q}(q^{\frac{1}{2}})$-subalgebra of $U_{q^{1/2}}(\widehat{\mathfrak{s l}}_2)$ generated by
\begin{eqnarray*}
	B_1=F_1+\varsigma_1 E_1 K_1^{-1};\hspace{0.2cm}
	 B_2=F_2+\varsigma_2 E_2 K_2^{-1} .
	\end{eqnarray*}
\end{definition}

\begin{corollary}\label{3.10}
	 The restriction $\Omega_{\varsigma}^{\imath}=\left.\Psi\right|_{U_{\varsigma}^{\imath}}: U_{\varsigma}^{\imath} \longrightarrow \mathcal{A}_q $ is an algebra homomorphism, where
\begin{eqnarray}
	&&\Omega_{\varsigma}^{\imath}\left(B_1\right)=w_1+q^{\frac{1}{2}} \varsigma_1 X^{(1,0,-1,1)};\nonumber \\
	 &&\Omega_{\boldsymbol{\varsigma}}^{\imath}\left(B_2\right)=X^{(0,1,0,-1)}-\frac{\varsigma_2}{\left(q^{\frac{1}{2}}-q^{-\frac{1}{2}}\right)^2}\left(q^{\frac{7}{2}} y_2^{-2} x_4+q y_2 x_0\right) . \label{(3.15)}
	\end{eqnarray}
\end{corollary}
\begin{proof}
	 Since $U_{\varsigma}^{\imath}$ is a subalgebra of $U_{q^{1/2}}(\widehat{\mathfrak{s l}}_2)$, the restriction of $\Psi$ to $U_{\varsigma}^{\imath}$ is an algebra homomorphism by Theorem \ref{3.6}. By (\ref{(3.6)}), (\ref{(w2)}) and (\ref{(ef)}), we get  (\ref{(3.15)}).
\end{proof}

\begin{proposition}
 We have $\mathbb{Q}(q^{\frac{1}{2}})\left\langle\operatorname{Im} \Omega_{\varsigma}^{\imath}, y_1^{ \pm 1}, y_2^{ \pm 1}\right\rangle=\mathcal{A}_q.$
\end{proposition}
\begin{proof}
	Set $\mathcal{S}=\mathbb{Q}(q^{\frac{1}{2}})\left\langle\operatorname{Im}
	\Omega_{\boldsymbol{\varsigma}}^{\imath}, y_1^{ \pm 1}, y_2^{ \pm 1}\right\rangle$. By Theorem \ref{3.8}, we have $\mathbb{Q}(q^{\frac{1}{2}})\left\langle\operatorname{Im} \Psi, y_1\right\rangle=\mathcal{A}_q$ It suffices to prove that $\Psi(K_1),\Psi(K_2),\Psi(E_1),\Psi(E_2),\Psi(F_1),\Psi(F_2) \in \mathcal{S}$.
	
Frist, we have $\Psi(K_1)=\Psi(K_2^{-1})=q^{-1} y_1 y_2^{-1};$ $ \Psi(K_2)=\Psi(K_1^{-1})=q^{-1} y_1^{-1} y_2;$ thus $\Psi(K_1), \Psi(K_2) \in \mathcal{S}.$

By Lemma \ref{3.3} and (\ref{(3.15)}), we obtain
\begin{eqnarray*}
	\Psi(K_1)\Omega_{\varsigma}^{\imath}\left(B_1\right) \Psi(K_2)-q^{-1} \Omega_{\varsigma}^{\imath}\left(B_1\right)&=& \varsigma_1\left(q-q^{-1}\right)\Psi(E_1)\Psi(K_2); \\
	\Psi(K_1) \Omega_{\varsigma}^{\imath}\left(B_2\right) \Psi(K_2)-q^{-1} \Omega_{\varsigma}^{\imath}\left(B_2\right)&=&\left(q-q^{-1}\right)\Psi(F_2) .
\end{eqnarray*}
Hence, $\Psi(E_1), \Psi(F_2) \in \mathcal{S}.$

By Definition \ref{4.1}, we have
	\begin{eqnarray*}
		\Omega_{\varsigma}^{\imath}\left(B_1\right)=\Psi(F_1)+ \varsigma_1\Psi(E_1)\Psi(K_1^{-1}) ;\hspace{0.2cm}
		\Omega_{\boldsymbol{\varsigma}}^{\imath}\left(B_2\right)=\Psi(F_2)+\varsigma_2\Psi(E_2)\Psi(K_2^{-1}) .
	\end{eqnarray*}
Therefore, $\Psi(E_2), \Psi(F_1) \in \mathcal{S}.$ This proof is complete.	
\end{proof}

 Let $\varsigma_1^{\diamond}=\varsigma_2^{\diamond}=-q^{-\frac{1}{2}}\left(q^{\frac{1}{2}}-q^{-\frac{1}{2}}\right)^2$, then the split $\imath$quantum group $U_{\varsigma^{\diamond}}^{\imath}$ is the $q^{1/2}$-Onsager algebra (see \cite[Remark~2.6]{LRZ25}). We denote it by $\mathcal{O}_{q^{1/2}}$.

\begin{corollary}
 The restriction $\Omega=\left.\Psi\right|_{\mathcal{O}_{q^{1/2}}}: \mathcal{O}_{q^{1/2}} \longrightarrow \mathcal{A}_q$ is an algebra homomorphism, where
\begin{eqnarray*}
	 &&\Omega\left(B_1\right)=w_1-\left(q^{\frac{1}{2}}-q^{-\frac{1}{2}}\right)^2 X^{(1,0,-1,1)}; \\
	 &&\Omega\left(B_2\right)=X^{(0,1,0,-1)}+q^3 y_2^{-2} x_4+q^{\frac{1}{2}} y_2 x_0 .
\end{eqnarray*}
\end{corollary}
\begin{proof}
	 This follows from Corollary \ref{3.10} and $\varsigma_1^{\diamond}=\varsigma_2^{\diamond}=-q^{-\frac{1}{2}}\left(q^{\frac{1}{2}}-q^{-\frac{1}{2}}\right)^2$.
\end{proof}

\section{The images of root vectors}

Recall the algebra homomorphism $\Psi: U_{q^{1/2}}(\widehat{\mathfrak{s l}}_2) \longrightarrow \mathcal{A}_q$ in Theorem \ref{3.6}. In this section, we give the images of the real and imaginary root vectors of $U_{q^{1/2}}(\widehat{\mathfrak{s l}}_2)$ under $\Psi$.

Let $\alpha_1, \alpha_2$ be the simple roots of type $A_1^{(1)}$ and  $\delta=\alpha_1+\alpha_2$. The set of positive real roots
 $$\bigtriangleup_{+}^{re}:=\{\alpha_1+n \delta,\ \alpha_2+n \delta \mid n \geq 0\},$$ and the set of  positive imaginary roots  $$\bigtriangleup_{+}^{im}:=\{n \delta, n \geq 1\}.$$
 The set of negative real roots is $\bigtriangleup_{-}^{re}=:-\bigtriangleup_{+}^{re}$, and the set of  negative imaginary roots is $\bigtriangleup_{-}^{im}:=-\bigtriangleup_{+}^{im}$.

 Following  Damiani \cite{Da93}, one can define the root vectors. Note that we use  $E_{\alpha}$ instead of $\widetilde{E}_{\alpha}$ in \cite{Da93} for any positive  root $\alpha$.

For the simple roots,	put $$E_{\alpha_1}=E_1,\ E_{\alpha_2}=E_2.$$
The root vector corresponding to the positive root $\delta$ is defined by
$$E_\delta=q^{-1} E_1 E_2-E_2 E_1.$$
For every $n \geq 1$,  the real root vectors corresponding to the positive
roots are defined recursively by
$$E_{\alpha_1+n \delta}=\frac{\left[E_{\alpha_1+(n-1) \delta}, E_\delta\right]}{q^{\frac{1}{2}}+q^{-\frac{1}{2}}};\hspace{0.2cm}  E_{\alpha_2+n \delta}=\frac{\left[E_\delta, E_{\alpha_2+(n-1) \delta}\right]}{q^{\frac{1}{2}}+q^{-\frac{1}{2}}};$$
and the positive imaginary root vectors corresponding to the positive imaginary roots are defined by
$$E_{n \delta}=q^{-1} E_{\alpha_1+(n-1) \delta} E_2-E_2 E_{\alpha_1+(n-1) \delta}.$$

The root vectors corresponding to the negative roots are defined as follows.

For the simple roots,	put $$F_{\alpha_1}=F_1;\hspace{0.2cm}  F_{\alpha_2}=F_2.$$
The root vector corresponding to the negative root $-\delta$ is defined by
$$F_\delta=q F_2 F_1-F_1 F_2.$$
For every $n \geq 1$,  the real root vectors corresponding to the negative
roots are defined recursively by
$$F_{\alpha_1+n \delta}=\frac{\left[F_\delta, F_{\alpha_1+(n-1) \delta}\right]}{q^{\frac{1}{2}}+q^{-\frac{1}{2}}};\hspace{0.2cm}  F_{\alpha_2+n \delta}=\frac{\left[F_{\alpha_2+(n-1) \delta}, F_\delta\right]}{q^{\frac{1}{2}}+q^{-\frac{1}{2}}},$$
and the  imaginary root vector corresponding to the negative imaginary roots  are defined by
$$F_{n \delta}=q F_2 F_{\alpha_1+(n-1) \delta}-F_{\alpha_1+(n-1) \delta} F_2.$$

\begin{proposition}
We have
\begin{flalign}
&\quad(1)\hspace{0.2cm} 	 \Psi\left(E_\delta\right)=\frac{q^{\frac{1}{2}}+q^{-\frac{1}{2}}}{q^{\frac{1}{2}}-q^{-\frac{1}{2}}} y_2^{-1} X_\delta;\label{(E'_delta)} && \\
&\quad(2)\hspace{0.2cm} 	 \Psi\left(F_\delta\right)=-\frac{q^{\frac{1}{2}}+q^{-\frac{1}{2}}}{q^{\frac{1}{2}}-q^{-\frac{1}{2}}} y_1^{-1} X_\delta \label{(F'_delta)}.&&
\end{flalign}
\end{proposition}
\begin{proof}
Applying the algebra homomorphism $\Psi$ to $E_\delta=q^{-1} E_1 E_2-E_2 E_1$  and $F_\delta=q F_2 F_1-F_1 F_2,$  we get
\begin{eqnarray*}
	\Psi\left(E_\delta\right)&=&q^{-1} x_1 w_2-w_2 x_1\mathop  = \limits^{(\ref{(w2x1-x1w2)})}\frac{q^{\frac{1}{2}}+q^{-\frac{1}{2}}}{q^{\frac{1}{2}}-q^{-\frac{1}{2}}} y_2^{-1} X_\delta;\\
	\Psi\left(F_\delta\right)&=&q X^{(0,1,0,-1)} w_1-w_1 X^{(0,1,0,-1)}\mathop  = \limits^{(\ref{(w1F2-F2w1)})}-\frac{q^{\frac{1}{2}}+q^{-\frac{1}{2}}}{q^{\frac{1}{2}}-q^{-\frac{1}{2}}} y_1^{-1} X_\delta .
	\end{eqnarray*}
The proof is complete.
\end{proof}
We need the following lemma which can be deduced from Theorem \ref{3.2}.
\begin{lemma}
In $\mathcal A_q$, we have
\begin{flalign}
&\quad(1)\hspace{0.2cm} 	x_{-m} X_\delta-q X_\delta x_{-m}  =(1-q) x_{-(m+1)},\hspace{0.2cm} \text{for $m\ge0$;} && \label{(x-mxd-xdx-m)}\\
&\quad(2)\hspace{0.2cm} 	X_\delta x_m-q x_m X_\delta  =(1-q) x_{m+1},\hspace{0.2cm} \text{for $m\ge1$}. &&\label{(xdxm-xmxd)}
\end{flalign}
\end{lemma}
\begin{theorem}
 For any $n \ge 1$, we have

\begin{flalign}
&\quad(1)\hspace{0.2cm} \Psi\left(E_{\alpha_1+n \delta}\right)=(-1)^n q^{\frac{n^2-4 n-1}{2}} y_1 y_2^{-n} x_{1-n};\label{(E'a1+nd)}&&\\
&\quad(2)\hspace{0.2cm}	\Psi\left(E_{\alpha_2+n \delta}\right)=\frac{(-1)^{n+1}}{\left(q^{\frac{1}{2}}-q^{-\frac{1}{2}}\right)^2}\left(q^{\frac{n^2+4 n+3}{2}} y_1^{-1} y_2^{-n-1} x_{n+4}+q^{\frac{n^2-4 n-1}{2}} y_1^{-1} y_2^{-n+3} x_n\right);\label{(E'a2+nd)}&&\\
&\quad(3)\hspace{0.2cm}	\Psi\left(E_{n \delta}\right)=  \frac{(-1)^{n+1}}{\left(q^{\frac{1}{2}}-q^{-\frac{1}{2}}\right)^2}\left(q^{\frac{n^2+2 n}{2}} y_1^{-1} y_2^{-n}\left(x_{n+3} x_1-x_1 x_{n+3}\right)+R(n)\right);\label{(E'nd)}&&
\end{flalign}

where $R(1)=0$ and $R(n)=q^{\frac{n^2-6 n+4}{2}} y_1^{-1} y_2^{-n+4}\left(x_{n-1} x_1-x_1 x_{n-1}\right)$ for any $n\geq 2$.
\end{theorem}
\begin{proof}
(1) We prove (\ref{(E'a1+nd)}) by induction on $n$. When $n=1$, by Theorem \ref{3.2} and applying the bar-involution, we obtain
$$\Psi\left(E_{\alpha_1+\delta}\right)=\frac{y_2^{-1}\left(x_1 X_\delta-X_\delta x_1\right)}{q^{\frac{1}{2}}-q^{-\frac{1}{2}}}=-y_2^{-1} y_1 x_0=-q^{-2} y_1 y_2^{-1} x_0 .$$
Suppose that this case holds for $n \le k$. When $n=k+1$, by (\ref{(3.8)}), and (\ref{(x-mxd-xdx-m)}), we have\begingroup
\interdisplaylinepenalty=0
\begin{eqnarray*}
	\Psi\left(E_{\alpha_1+(k+1) \delta}\right)&=&\frac{{\left[ {\Psi \left( {{E_{{\alpha _1} + k\delta }}} \right),\Psi \left( {{E_\delta }} \right)} \right]}}{{{q^{\frac{1}{2}}} + {q^{ - \frac{1}{2}}}}}  \mathop  = \limits^{ (\ref{(E'_delta)})}\frac{\left[\Psi\left(E_{\alpha_1+k \delta}\right), y_2^{-1} X_\delta\right]}{q^{\frac{1}{2}}-q^{-\frac{1}{2}}} \\
	& =&\frac{1}{q^{\frac{1}{2}}-q^{-\frac{1}{2}}}\left(\Psi\left(E_{\alpha_1+k \delta}\right) y_2^{-1} X_\delta-y_2^{-1} X_\delta \Psi\left(E_{\alpha_1+k \delta}\right)\right) \\
	& =&\frac{(-1)^k q^{\frac{k^2-4 k-1}{2}}}{q^{\frac{1}{2}}-q^{-\frac{1}{2}}}\left(y_1 y_2^{-k} x_{1-k} y_2^{-1} X_\delta-y_2^{-1} X_\delta y_1 y_2^{-k} x_{1-k}\right) \\
	& =&\frac{(-1)^k q^{\frac{k^2-4 k-1}{2}}}{q^{\frac{1}{2}}-q^{-\frac{1}{2}}} y_1 y_2^{-k-1}\left(q^{k-2} x_{1-k} X_\delta-q^{k-1} X_\delta x_{1-k}\right) \\
	& =&\frac{(-1)^k q^{\frac{k^2-2 k-5}{2}}}{q^{\frac{1}{2}}-q^{-\frac{1}{2}}} y_1 y_2^{-k-1}\left(x_{1-k} X_\delta-q X_\delta x_{1-k}\right) \\
	& \mathop  = \limits^{(\ref{(x-mxd-xdx-m)})}&\frac{(-1)^k q^{\frac{k^2-2 k-5}{2}}}{q^{\frac{1}{2}}-q^{-\frac{1}{2}}} y_1 y_2^{-k-1}(1-q) x_{-k}=(-1)^{k+1} q^{\frac{k^2-2 k-4}{2}} y_1 y_2^{-k-1} x_{-k} .
\end{eqnarray*}

(2) By Theorem \ref{3.2} and the bar-involution, we have $X_\delta x_0-x_0 X_\delta=\left(q^{-\frac{1}{2}}-\right.$ $\left.q^{\frac{1}{2}}\right) x_1 y_2$.

We prove (\ref{(E'a2+nd)}) by induction on $n$. When $n=1$, by (\ref{(w2)}) and  (\ref{(xdxm-xmxd)}), we obtain
\begin{eqnarray*}
	\Psi\left(E_{\alpha_2+\delta}\right)&= & \frac{\left[y_2^{-1} X_\delta, w_2\right]}{q^{\frac{1}{2}}-q^{-\frac{1}{2}}}=-\frac{1}{\left(q^{\frac{1}{2}}-q^{-\frac{1}{2}}\right)^3}\left(q^{\frac{7}{2}} y_1^{-1} y_2^{-2}\left(X_\delta x_4-q x_4 X_\delta\right)\right. \\
	&& \left.+q^{-2} y_1^{-1} y_2\left(X_\delta x_0-x_0 X_\delta\right)\right) \\
	& \mathop  = \limits^{(\ref{(xdxm-xmxd)})} & -\frac{1}{\left(q^{\frac{1}{2}}-q^{-\frac{1}{2}}\right)^3}\left((1-q) q^{\frac{7}{2}} y_1^{-1} y_2^{-2} x_5-\left(q^{\frac{1}{2}}-q^{-\frac{1}{2}}\right) q^{-2} y_1^{-1} y_2 x_1 y_2\right) \\
&	= & \frac{q^4 y_1^{-1} y_2^{-2} x_5+q^{-2} y_1^{-1} y_2^2 x_1}{\left(q^{\frac{1}{2}}-q^{-\frac{1}{2}}\right)^2} .
\end{eqnarray*}

Suppose that this case holds for $n \le k$. When $n=k+1$, by (\ref{(3.2)}) and (\ref{(xdxm-xmxd)}), we have
\begin{eqnarray*}
	&& \Psi\left(E_{\alpha_2+(k+1) \delta}\right)=\frac{{\left[ {\Psi \left( {{E_\delta }} \right),\Psi \left( {{E_{{\alpha _2} + k\delta }}} \right)} \right]}}{{{q^{\frac{1}{2}}} + {q^{ - \frac{1}{2}}}}} \\
&	 \mathop  = \limits^{ (\ref{(E'_delta)})} & \frac{\left[y_2^{-1} X_\delta, \Psi\left(E_{\alpha_2+k \delta}\right)\right]}{q^{\frac{1}{2}}-q^{-\frac{1}{2}}}=\frac{1}{q^{\frac{1}{2}}-q^{-\frac{1}{2}}}\left(y_2^{-1} X_\delta \Psi\left(E_{\alpha_2+k \delta}\right)-\Psi\left(E_{\alpha_2+k \delta}\right) y_2^{-1} X_\delta\right) \\
&	= & \frac{(-1)^{k+1}}{\left(q^{\frac{1}{2}}-q^{-\frac{1}{2}}\right)^3}\left(y_2^{-1} X_\delta\left(q^{\frac{k^2+4 k+3}{2}} y_1^{-1} y_2^{-k-1} x_{k+4}+q^{\frac{k^2-4 k-1}{2}} y_1^{-1} y_2^{-k+3} x_k\right)\right. \\
&	& \left.-\left(q^{\frac{k^2+4 k+3}{2}} y_1^{-1} y_2^{-k-1} x_{k+4}+q^{\frac{k^2-4 k-1}{2}} y_1^{-1} y_2^{-k+3} x_k\right) y_2^{-1} X_\delta\right) \\
&	= & \frac{(-1)^{k+1}}{\left(q^{\frac{1}{2}}-q^{-\frac{1}{2}}\right)^3}\left(q^{\frac{k^2+6 k+7}{2}} y_1^{-1} y_2^{-k-2}\left(X_\delta x_{k+4}-q x_{k+4} X_\delta\right)\right. \\
&	& \left.+q^{\frac{k^2-2 k-5}{2}} y_1^{-1} y_2^{-k+2}\left(X_\delta x_k-q x_k X_\delta\right)\right) \\
&	 \mathop  = \limits^{(\ref{(xdxm-xmxd)})}& \frac{(-1)^{k+2}}{\left(q^{\frac{1}{2}}-q^{-\frac{1}{2}}\right)^2}\left(q^{\frac{k^2+6 k+8}{2}} y_1^{-1} y_2^{-k-2} x_{k+5}+q^{\frac{k^2-2 k-4}{2}} y_1^{-1} y_2^{-k+2} x_{k+1}\right) .
\end{eqnarray*}\endgroup

(3) By \cite{Ter19}, for any $r, s \ge 0$, we have
$$q^{-1} E_{\alpha_1+r \delta} E_{\alpha_2+s \delta}-E_{\alpha_2+s\delta} E_{\alpha_1+r \delta}=E_{(r+s+1) \delta} .$$

When $r=0$ and $s=n-1$, by $\Psi$, we obtain
$$\Psi\left(E_{n \delta}\right)=q^{-1} x_1 \Psi\left(E_{\alpha_2+(n-1) \delta}\right)-\Psi\left(E_{\alpha_2+(n-1) \delta}\right) x_1 .$$

For $n \ge 2$, substitute (\ref{(E'a2+nd)}). Note that $x_1 y_1^{-1} y_2^b=q y_1^{-1} y_2^b x_1$, we obtain
\begingroup
\interdisplaylinepenalty=0
\begin{eqnarray*}\Psi\left(E_{n \delta}\right)&=&q^{-1} x_1 \Psi\left(E_{\alpha_2+(n-1) \delta}\right)-\Psi\left(E_{\alpha_2+(n-1) \delta}\right) x_1\\
&	= & \frac{(-1)^n}{\left(q^{\frac{1}{2}}-q^{-\frac{1}{2}}\right)^2}\left(q^{\frac{n^2+2 n}{2}-1} x_1 y_1^{-1} y_2^{-n} x_{n+3}+q^{\frac{n^2-6 n+4}{2}-1} x_1 y_1^{-1} y_2^{-n+4} x_{n-1}\right. \\
&	& \left.-q^{\frac{n^2+2 n}{2}} y_1^{-1} y_2^{-n} x_{n+3} x_1-q^{\frac{n^2-6 n+4}{2}} y_1^{-1} y_2^{-n+4} x_{n-1} x_1\right) \\
&	= & \frac{(-1)^n}{\left(q^{\frac{1}{2}}-q^{-\frac{1}{2}}\right)^2}\left(q^{\frac{n^2+2 n}{2}} y_1^{-1} y_2^{-n} x_1 x_{n+3}+q^{\frac{n^2-6 n+4}{2}} y_1^{-1} y_2^{-n+4} x_1 x_{n-1}\right. \\
&	& \left.-q^{\frac{n^2+2 n}{2}} y_1^{-1} y_2^{-n} x_{n+3} x_1-q^{\frac{n^2-6 n+4}{2}} y_1^{-1} y_2^{-n+4} x_{n-1} x_1\right) \\
&	= & \frac{(-1)^n}{\left(q^{\frac{1}{2}}-q^{-\frac{1}{2}}\right)^2}\left(q^{\frac{n^2+2 n}{2}} y_1^{-1} y_2^{-n}\left(x_1 x_{n+3}-x_{n+3} x_1\right)\right. \\
&	& \left.+q^{\frac{n^2-6 n+4}{2}} y_1^{-1} y_2^{-n+4}\left(x_1 x_{n-1}-x_{n-1} x_1\right)\right) \\
&	= & \frac{(-1)^{n+1}}{\left(q^{\frac{1}{2}}-q^{-\frac{1}{2}}\right)^2}\left(q^{\frac{n^2+2 n}{2}} y_1^{-1} y_2^{-n}\left(x_{n+3} x_1-x_1 x_{n+3}\right)\right. \\
&	& \left.+q^{\frac{n^2-6 n+4}{2}} y_1^{-1} y_2^{-n+4}\left(x_{n-1} x_1-x_1 x_{n-1}\right)\right) .
\end{eqnarray*}
\endgroup

For $n=1$, note that $x_4 x_1-x_1 x_4 =\left(q^{\frac{3}{2}}-q^{-\frac{1}{2}}\right) y_1 X_\delta$. By (\ref{(E'_delta)}), we obtain
\begin{eqnarray*}
	\Psi\left(E_\delta\right)&= &  \frac{{q\left( {1 - {q^{ - 2}}} \right)}}{{{{\left( {{q^{\frac{1}{2}}} - {q^{ - \frac{1}{2}}}} \right)}^2}}}y_2^{ - 1}{X_\delta } = \frac{{{q^{ - \frac{1}{2}}}}}{{{{\left( {{q^{\frac{1}{2}}} - {q^{ - \frac{1}{2}}}} \right)}^2}}}y_2^{ - 1}y_1^{ - 1}\left( {{x_4}{x_1} - {x_1}{x_4}} \right)\\
&	= &\frac{{{q^{\frac{3}{2}}}}}{{{{\left( {{q^{\frac{1}{2}}} - {q^{ - \frac{1}{2}}}} \right)}^2}}}y_1^{ - 1}y_2^{ - 1}\left( {{x_4}{x_1} - {x_1}{x_4}} \right).
\end{eqnarray*}
The proof is complete.
\end{proof}

Note that, we can rewrite $\Psi \left( {{E_{n\delta }}} \right)$ in terms of $S_n$ $(n\geq 0)$. To do this, we need the following lemma.
\begin{lemma}\cite{DXC23}\label{5.4}
	For any $n \ge 0$, we have
\begin{eqnarray*}
		{y_1}{S_n} = {q^{ - n}}{S_n}{y_1};\hspace{0.2cm}
		{X^{(0,0,1,1)}}{S_n} = {S_n}{X^{(0,0,1,1)}}.
\end{eqnarray*}
\end{lemma}

	For any $j\in \mathbb{N}$, set
\[\varepsilon_j=
\begin{cases}
	0, & j=0,1;\\
	1, & j\ge 2.
\end{cases}\]
\begin{corollary}\label{5.5}
	For any $n\ge 1$, we have
	\begin{eqnarray*}
	\Psi \left( {{E_{n\delta }}} \right) = \frac{{{{( - 1)}^{n + 1}}{q^{\frac{{{n^2} - 2n - 1}}{2}}}}}{{{{\left( {{q^{\frac{1}{2}}} - {q^{ - \frac{1}{2}}}} \right)}^2}}}y_2^{ - n}\sum\limits_{j = 0}^{\left\lfloor {\frac{n}{2}} \right\rfloor } {\left( {{q^{n - j + 1}} - {q^j}} \right)\left( {1 + {\varepsilon _j}{q^{ - 4n - 4}}y_1^{ - 2}y_2^2} \right){X^{(0,0,j,j)}}{S_{n - 2j}}} .
	\end{eqnarray*}
\end{corollary}
\begin{proof}
	By Theorem \ref{3.1} (2) and  the bar-involution, for any $m\ge 1$ and $n\ge 0$, we have\begingroup
	\interdisplaylinepenalty=0
\begin{eqnarray}\label{(5.9)}
	&&	x_{m+n+2}x_m-x_mx_{m+n+2}\nonumber\\
	&=&\left(q^{\frac12}-q^{-n-\frac12}\right)
		S_nX^{(0,0,m,m-1)}+x_{m+n+1}x_{m+1}-x_{m+1}x_{m+n+1}.
\end{eqnarray}

By repeatedly applying (\ref{(5.9)}) for $0\le j\le \left\lfloor\frac{n}{2}\right\rfloor$, we obtain
\begin{eqnarray*}
&&	{x_{n + 3}}{x_1} - {x_1}{x_{n + 3}} \\
&=& \left( {{q^{\frac{1}{2}}} - {q^{ - n - \frac{1}{2}}}} \right){S_n}{y_1} + {x_{n + 2}}{x_2} - {x_2}{x_{n + 2}}\\
&\mathop  = \limits^{(\ref{(5.9)})} &\sum\limits_{j = 0}^{ \left\lfloor {\frac{{n }}{2}} \right\rfloor } {\left( {{q^{\frac{1}{2}}} - {q^{ - n + 2j - \frac{1}{2}}}} \right){S_{n - 2j}}{X^{(0,0,j + 1,j)}}}\\
&&  + {x_{n +2- \left\lfloor {\frac{{n}}{2}} \right\rfloor }}{x_{2+ \left\lfloor {\frac{{n }}{2}} \right\rfloor }} - {x_{2 + \left\lfloor {\frac{{n }}{2}} \right\rfloor }}{x_{n+2 - \left\lfloor {\frac{{n }}{2}} \right\rfloor }}\\
&	=& \sum\limits_{j = 0}^{\left\lfloor {\frac{n}{2}} \right\rfloor } {\left( {{q^{\frac{1}{2}}} - {q^{ - n + 2j - \frac{1}{2}}}} \right){q^j}{S_{n - 2j}}{y_1}{X^{(0,0,j,j)}}}  + {x_{\left\lceil {\frac{{n + 4}}{2}} \right\rceil }}{x_{\left\lfloor {\frac{{n + 4}}{2}} \right\rfloor }} - {x_{\left\lfloor {\frac{{n + 4}}{2}} \right\rfloor }}{x_{\left\lceil {\frac{{n + 4}}{2}} \right\rceil }}\\
&	= &{y_1}\sum\limits_{j = 0}^{\left\lfloor {\frac{n}{2}} \right\rfloor } {{q^{n - j}}\left( {{q^{\frac{1}{2}}} - {q^{ - n + 2j - \frac{1}{2}}}} \right){X^{(0,0,j,j)}}{S_{n - 2j}}} \\
	&=& {q^{ - \frac{1}{2}}}{y_1}\sum\limits_{j = 0}^{\left\lfloor {\frac{n}{2}} \right\rfloor } {\left( {{q^{n - j + 1}} - {q^j}} \right){X^{(0,0,j,j)}}{S_{n - 2j}}} .
\end{eqnarray*}
	Similarly, we obtain
\[{x_{n - 1}}{x_1} - {x_1}{x_{n - 1}} = {q^{ - \frac{1}{2}}}{y_1}\sum\limits_{j = 0}^{\left\lfloor {\frac{{n - 4}}{2}} \right\rfloor } {\left( {{q^{n - j - 3}} - {q^j}} \right){X^{(0,0,j,j)}}{S_{n - 4 - 2j}}}. \]
	
	By (\ref{(E'nd)}), Theorem \ref{3.1} (2) and the bar-involution, for any $n\geq 4$, we obtain
\begin{eqnarray*}
	\Psi \left( {{E_{n\delta }}} \right)& =& \frac{{{{( - 1)}^{n + 1}}}}{{{{\left( {{q^{\frac{1}{2}}} - {q^{ - \frac{1}{2}}}} \right)}^2}}}\left( {{q^{\frac{{{n^2} - 2n - 1}}{2}}}y_2^{ - n}\sum\limits_{j = 0}^{\left\lfloor {\frac{n}{2}} \right\rfloor } {\left( {{q^{n - j + 1}} - {q^j}} \right){X^{(0,0,j,j)}}{S_{n - 2j}}} } \right.\\
&&	\left. { + {q^{\frac{{{n^2} - 10n + 19}}{2}}}y_2^{ - n + 4}\sum\limits_{j = 0}^{\left\lfloor {\frac{{n - 4}}{2}} \right\rfloor } {\left( {{q^{n - j - 3}} - {q^j}} \right){X^{(0,0,j,j)}}{S_{n - 4 - 2j}}} } \right)\\
&	=& \frac{{{{( - 1)}^{n + 1}}{q^{\frac{{{n^2} - 2n - 1}}{2}}}}}{{{{\left( {{q^{\frac{1}{2}}} - {q^{ - \frac{1}{2}}}} \right)}^2}}}y_2^{ - n}\left( {\sum\limits_{j = 0}^{\left\lfloor {\frac{n}{2}} \right\rfloor } {\left( {{q^{n - j + 1}} - {q^j}} \right){X^{(0,0,j,j)}}{S_{n - 2j}}} } \right.\\
&&	\left. { + {q^{ - 4n + 10}}y_2^4\sum\limits_{t = 2}^{\left\lfloor {\frac{n}{2}} \right\rfloor } {\left( {{q^{n - t - 1}} - {q^{t - 2}}} \right){X^{(0,0,t - 2,t - 2)}}{S_{n - 2t}}} } \right)\\
&	= &\frac{{{{( - 1)}^{n + 1}}{q^{\frac{{{n^2} - 2n - 1}}{2}}}}}{{{{\left( {{q^{\frac{1}{2}}} - {q^{ - \frac{1}{2}}}} \right)}^2}}}y_2^{ - n}\left( {\sum\limits_{j = 0}^{\left\lfloor {\frac{n}{2}} \right\rfloor } {\left( {{q^{n - j + 1}} - {q^j}} \right){X^{(0,0,j,j)}}{S_{n - 2j}}} } \right.\\
&&	\left. { + {q^{ - 4n + 10}}\sum\limits_{t = 2}^{\left\lfloor {\frac{n}{2}} \right\rfloor } {{q^{  - 12}}\left( {{q^{n - t - 1}} - {q^{t - 2}}} \right)y_1^{ - 2}y_2^2{X^{(0,0,t,t)}}{S_{n - 2t}}} } \right)\\
&	=&\frac{{{{( - 1)}^{n + 1}}{q^{\frac{{{n^2} - 2n - 1}}{2}}}}}{{{{\left( {{q^{\frac{1}{2}}} - {q^{ - \frac{1}{2}}}} \right)}^2}}}y_2^{ - n}\left( {\sum\limits_{j = 0}^{\left\lfloor {\frac{n}{2}} \right\rfloor } {\left( {{q^{n - j + 1}} - {q^j}} \right){X^{(0,0,j,j)}}{S_{n - 2j}}} } \right.\\
&&	\left. { + \sum\limits_{t = 2}^{\left\lfloor {\frac{n}{2}} \right\rfloor } {{q^{ - 4n  - 4}}\left( {{q^{n - t + 1}} - {q^t}} \right)y_1^{ - 2}y_2^2{X^{(0,0,t,t)}}{S_{n - 2t}}} } \right)\\
&	=& \frac{{{{( - 1)}^{n + 1}}{q^{\frac{{{n^2} - 2n - 1}}{2}}}}}{{{{\left( {{q^{\frac{1}{2}}} - {q^{ - \frac{1}{2}}}} \right)}^2}}}y_2^{ - n}\sum\limits_{j = 0}^{\left\lfloor {\frac{n}{2}} \right\rfloor } {\left( {{q^{n - j + 1}} - {q^j}} \right)\left( {1 + \varepsilon_j{q^{ - 4n  - 4}}y_1^{ - 2}y_2^2} \right){X^{(0,0,j,j)}}{S_{n - 2j}}} .
\end{eqnarray*}\endgroup
	For $n=1,2,3$, note that $R(n)=0$. Thus, we have
	\[\Psi \left( {{E_{n\delta }}} \right) = {\rm{ }}\frac{{{{( - 1)}^{n + 1}}{q^{\frac{{{n^2} - 2n - 1}}{2}}}}}{{{{\left( {{q^{\frac{1}{2}}} - {q^{ - \frac{1}{2}}}} \right)}^2}}}y_2^{ - n}\sum\limits_{j = 0}^{\left\lfloor {\frac{n}{2}} \right\rfloor } {\left( {{q^{n - j + 1}} - {q^j}} \right){X^{(0,0,j,j)}}{S_{n - 2j}}} .\]
	This completes the proof.
\end{proof}

\begin{theorem}\label{5.6}
 (1) For any $n \ge 1$, we have
\begin{equation}\label{(F'a2+nd)}
	\Psi\left(F_{\alpha_2+n \delta}\right)=(-1)^n q^{\frac{-n^2+4 n+1}{2}} y_1^{-n} y_2^{-1} x_{n+2} .
\end{equation}

(2) For any $n \ge 3$, we have
\begin{equation}\label{(F'a1+nd)}
	\Psi\left(F_{\alpha_1+n \delta}\right)=\frac{(-1)^{n+1}}{\left(q^{\frac{1}{2}}-q^{-\frac{1}{2}}\right)^2}\left(q^{-\frac{n^2+4 n+3}{2}} y_1^{-n-1} y_2 x_{-(n+1)}+q^{\frac{-n^2+4 n+1}{2}} y_1^{-n+3} y_2 x_{3-n}\right) .
\end{equation}

(3) For any $n \ge 1$, we have
\begin{eqnarray}\label{(F'nd)}
	\Psi\left(F_{n \delta}\right)&= & \frac{(-1)^{n+1}}{\left(q^{\frac{1}{2}}-q^{-\frac{1}{2}}\right)^2}\left(q^{-\frac{n^2+1}{2}} y_1^{-n}\left(x_{-n} x_2-q^{n+1} x_2 x_{-n}\right)+R'(n)\right),
\end{eqnarray}	
where $R'(1)=R'(2)=R'(3)=0$ and $R'(n)=q^{\frac{3-(n-4)^2}{2}} y_1^{-n+4}\left(x_{4-n} x_2-q^{n-3} x_2 x_{4-n}\right)$ for any $n\ge 4$.
\end{theorem}
\begin{proof}
(1) We prove (\ref{(F'a2+nd)}) by induction on $n$. When $n=1$, by (\ref{(3.6)}), (\ref{(ef)}) and (\ref{(xdxm-xmxd)}), we obtain\begingroup
\interdisplaylinepenalty=0
\begin{eqnarray*}
	&&\Psi\left(F_{\alpha_2+\delta}\right)  =\frac{\left[y_1^{-1} X_\delta, q^{\frac{1}{2}} y_2^{-1} x_2\right]}{q^{\frac{1}{2}}-q^{-\frac{1}{2}}}=\frac{q^{\frac{1}{2}}}{q^{\frac{1}{2}}-q^{-\frac{1}{2}}}\left(y_1^{-1} X_\delta y_2^{-1} x_2-y_2^{-1} x_2 y_1^{-1} X_\delta\right) \\
	& =&\frac{q^{\frac{3}{2}}}{q^{\frac{1}{2}}-q^{-\frac{1}{2}}} y_1^{-1} y_2^{-1}\left(X_\delta x_2-q x_2 X_\delta\right)=\frac{q^{\frac{3}{2}}(1-q)}{q^{\frac{1}{2}}-q^{-\frac{1}{2}}} y_1^{-1} y_2^{-1} x_3=-q^2 y_1^{-1} y_2^{-1} x_3 .
\end{eqnarray*}
Suppose that the case holds for $1 \le n \le k$. When $n=k+1$, by (\ref{(3.8)}), (\ref{(3.2)}) and (\ref{(xdxm-xmxd)}), we have
\begin{eqnarray*}
	&&\Psi\left(F_{\alpha_2+(k+1) \delta}\right)  =\frac{{\left[ {\Psi \left( {{F_{{\alpha _2} + k\delta }}} \right),\Psi \left( {{F_\delta }} \right)} \right]}}{{{q^{\frac{1}{2}}} + {q^{ - \frac{1}{2}}}}}\mathop  = \limits^{(\ref{(F'_delta)})}\frac{\left[y_1^{-1} X_\delta, \Psi\left(F_{\alpha_2+k \delta}\right)\right]}{q^{\frac{1}{2}}-q^{-\frac{1}{2}}} \\
	& =&\frac{1}{q^{\frac{1}{2}}-q^{-\frac{1}{2}}}\left(y_1^{-1} X_\delta \Psi\left(F_{\alpha_2+k \delta}\right)-\Psi\left(F_{\alpha_2+k \delta}\right) y_1^{-1} X_\delta\right) \\
	& =&\frac{(-1)^k q^{\frac{-k^2+4 k+1}{2}}}{q^{\frac{1}{2}}-q^{-\frac{1}{2}}}\left(y_1^{-1} X_\delta y_1^{-k} y_2^{-1} x_{k+2}-y_1^{-k} y_2^{-1} x_{k+2} y_1^{-1} X_\delta\right) \\
	& =&\frac{(-1)^k q^{\frac{-k^2+4 k+1}{2}}}{q^{\frac{1}{2}}-q^{-\frac{1}{2}}} y_1^{-k-1} y_2^{-1}\left(q^{1-k} X_\delta x_{k+2}-q^{2-k} x_{k+2} X_\delta\right) \\
	& 	=&\frac{(-1)^k q^{\frac{-k^2+2 k+3}{2}}}{q^{\frac{1}{2}}-q^{-\frac{1}{2}}} y_1^{-k-1} y_2^{-1}\left(X_\delta x_{k+2}-q x_{k+2} X_\delta\right) \\
	& \mathop  = \limits^{(\ref{(xdxm-xmxd)})}& \frac{(-1)^k q^{\frac{-k^2+2 k+3}{2}}}{q^{\frac{1}{2}}-q^{-\frac{1}{2}}}(1-q) y_1^{-k-1} y_2^{-1} x_{k+3} \\
	& =&(-1)^{k+1} q^{\frac{-k^2+2 k+4}{2}} y_1^{-k-1} y_2^{-1} x_{k+3} .
\end{eqnarray*}

(2) By Theorem \ref{3.2} and (\ref{(xdxm-xmxd)}), we have
\begin{eqnarray*}
	x_3 X_\delta-X_\delta x_3  =x_3 X_\delta-q x_3 X_\delta-(1-q) x_4
	 =(1-q)\left(x_3 X_\delta-x_4\right)=(1-q) y_1 y_2 x_2 .
\end{eqnarray*}
Thus, by (\ref{(w1)}) and (\ref{(x-mxd-xdx-m)}), we obtain
\begin{eqnarray*}
	&&\Psi\left(F_{\alpha_1+\delta}\right)  =\frac{\left[w_1, y_1^{-1} X_\delta\right]}{q^{\frac{1}{2}}-q^{-\frac{1}{2}}} \\
	& =&-\frac{1}{\left(q^{\frac{1}{2}}-q^{-\frac{1}{2}}\right)^3}\left(q^{-\frac{9}{2}} y_1^{-2} y_2\left(x_{-1} X_\delta-q X_\delta x_{-1}\right)+q^2 y_1^{-1} y_2^{-1}\left(x_3 X_\delta-X_\delta x_3\right)\right) \\
	& =&-\frac{1-q}{\left(q^{\frac{1}{2}}-q^{-\frac{1}{2}}\right)^3}\left(q^{-\frac{9}{2}} y_1^{-2} y_2 x_{-2}+q^2 y_1^{-1} y_2^{-1} y_1 y_2 x_2\right) \\
	& =&\frac{q^{-4} y_1^{-2} y_2 x_{-2}+q^{\frac{1}{2}} x_2}{\left(q^{\frac{1}{2}}-q^{-\frac{1}{2}}\right)^2} .
\end{eqnarray*}
By Theorem \ref{3.2} and the bar-involution and (\ref{(x-mxd-xdx-m)}), we get
\begin{eqnarray*}
	\Psi\left(F_{\alpha_1+2 \delta}\right)=\frac{{\left[ {\Psi \left( {{F_\delta }} \right),\Psi \left( {{F_{{\alpha _1} + \delta }}} \right)} \right]}}{{{q^{\frac{1}{2}}} + {q^{ - \frac{1}{2}}}}} \mathop  = \limits^{(\ref{(F'_delta)})}\frac{\left[\Psi\left(F_{\alpha_1+\delta}\right), y_1^{-1} X_\delta\right]}{q^{\frac{1}{2}}-q^{-\frac{1}{2}}}=-\frac{q^{-\frac{15}{2}} y_1^{-3} y_2 x_{-3}+q y_2 x_1}{\left(q^{\frac{1}{2}}-q^{-\frac{1}{2}}\right)^2}; \\
	 \Psi\left(F_{\alpha_1+3 \delta}\right)=\frac{{\left[ {\Psi \left( {{F_\delta }} \right),\Psi \left( {{F_{{\alpha _1} +2 \delta }}} \right)} \right]}}{{{q^{\frac{1}{2}}} + {q^{ - \frac{1}{2}}}}} \mathop  = \limits^{(\ref{(F'_delta)})}\frac{\left[\Psi\left(F_{\alpha_1+2 \delta}\right), y_1^{-1} X_\delta\right]}{q^{\frac{1}{2}}-q^{-\frac{1}{2}}}=\frac{q^{-12} y_1^{-4} y_2 x_{-4}+q^2 y_2 x_0}{\left(q^{\frac{1}{2}}-q^{-\frac{1}{2}}\right)^2} .
\end{eqnarray*}
Thus (\ref{(F'a1+nd)}) holds for $n=3$. Suppose that this case holds for $3 \le n \le k$. When $n=k+1$, by (\ref{(3.1)}) and (\ref{(x-mxd-xdx-m)}), we have
\begin{eqnarray*}
	&& \Psi\left(F_{\alpha_1+(k+1) \delta}\right)=\frac{{\left[ {\Psi \left( {{F_\delta }} \right),\Psi \left( {{F_{{\alpha _1} + k\delta }}} \right)} \right]}}{{{q^{\frac{1}{2}}} + {q^{ - \frac{1}{2}}}}}
 \mathop  = \limits^{(\ref{(F'_delta)})} \frac{\left[\Psi\left(F_{\alpha_1+k \delta}\right), y_1^{-1} X_\delta\right]}{q^{\frac{1}{2}}-q^{-\frac{1}{2}}} \\
&	= & \frac{1}{q^{\frac{1}{2}}-q^{-\frac{1}{2}}}\left(\Psi\left(F_{\alpha_1+k \delta}\right) y_1^{-1} X_\delta-y_1^{-1} X_\delta \Psi\left(F_{\alpha_1+k \delta}\right)\right) \\
&	= & \frac{(-1)^{k+1}}{\left(q^{\frac{1}{2}}-q^{-\frac{1}{2}}\right)^3}\left(q^{-\frac{k^2+4 k+3}{2}}\left(y_1^{-k-1} y_2 x_{-(k+1)} y_1^{-1} X_\delta-y_1^{-1} X_\delta y_1^{-k-1} y_2 x_{-(k+1)}\right)\right. \\
	&& \left.+q^{\frac{-k^2+4 k+1}{2}}\left(y_1^{-k+3} y_2 x_{3-k} y_1^{-1} X_\delta-y_1^{-1} X_\delta y_1^{-k+3} y_2 x_{3-k}\right)\right) \\
&	= & \frac{(-1)^{k+1}}{\left(q^{\frac{1}{2}}-q^{-\frac{1}{2}}\right)^3}\left(q^{-\frac{k^2+4 k+3}{2}} y_1^{-k-2} y_2\left(q^{-k-3} x_{-(k+1)} X_\delta-q^{-k-2} X_\delta x_{-(k+1)}\right)\right. \\
&	& \left.+q^{\frac{-k^2+4 k+1}{2}} y_1^{-k+2} y_2\left(q^{-k+1} x_{3-k} X_\delta-q^{-k+2} X_\delta x_{3-k}\right)\right) \\
&	= & \frac{(-1)^{k+1}}{\left(q^{\frac{1}{2}}-q^{-\frac{1}{2}}\right)^3}\left(q^{-\frac{k^2+6 k+9}{2}} y_1^{-k-2} y_2\left(x_{-(k+1)} X_\delta-q X_\delta x_{-(k+1)}\right)\right. \\
&	& \left.+q^{\frac{-k^2+2 k+3}{2}} y_1^{-k+2} y_2\left(x_{3-k} X_\delta-q X_\delta x_{3-k}\right)\right)\\
	&  \mathop  = \limits^{(\ref{(x-mxd-xdx-m)})}&\frac{(-1)^{k+1}(1-q)}{\left(q^{\frac{1}{2}}-q^{-\frac{1}{2}}\right)^3}\left(q^{-\frac{k^2+6 k+9}{2}} y_1^{-k-2} y_2 x_{-(k+2)}+q^{\frac{-k^2+2 k+3}{2}} y_1^{-k+2} y_2 x_{2-k}\right) \\
& =&\frac{(-1)^{k+2}}{\left(q^{\frac{1}{2}}-q^{-\frac{1}{2}}\right)^2}\left(q^{-\frac{k^2+6 k+8}{2}} y_1^{-k-2} y_2 x_{-(k+2)}+q^{\frac{-k^2+2 k+4}{2}} y_1^{-k+2} y_2 x_{2-k}\right) .
\end{eqnarray*}
\endgroup

(3) First, we have
$$\Psi\left(F_{n \delta}\right)=q^{\frac{3}{2}} y_2^{-1} x_2 \Psi\left(F_{\alpha_1+(n-1) \delta}\right)-q^{\frac{1}{2}} \Psi\left(F_{\alpha_1+(n-1) \delta}\right) y_2^{-1} x_2 .$$
For any $n \ge 4$, by (\ref{(F'a1+nd)}), we have
$$\Psi\left(F_{\alpha_1+(n-1) \delta}\right)=\frac{(-1)^n}{\left(q^{\frac{1}{2}}-q^{-\frac{1}{2}}\right)^2}\left(q^{-\frac{n^2+2 n}{2}} y_1^{-n} y_2 x_{-n}+q^{\frac{-n^2+6 n-4}{2}} y_1^{-n+4} y_2 x_{4-n}\right) .$$

Note that, by (\ref{(ef)}) and (\ref{(3.1)}), for any $m \ge 0$, we have
\begin{eqnarray*}
	& y_2^{-1} x_2 y_1^{-m} y_2 x_{-m}=q^{2 m} y_1^{-m} y_2^{-1} x_2 y_2 x_{-m}=q^{2 m-1} y_1^{-m} x_2 x_{-m}; \\
	& y_1^{-m} y_2 x_{-m} y_2^{-1} x_2=q^{m-1} y_1^{-m} y_2 y_2^{-1} x_{-m} x_2=q^{m-1} y_1^{-m} x_{-m} x_2 .
\end{eqnarray*}
Thus, we obtain
\begingroup
\interdisplaylinepenalty=0
\begin{eqnarray*}
	&&\Psi\left(F_{n \delta}\right)=  q^{\frac{3}{2}} y_2^{-1} x_2 \Psi\left(F_{\alpha_1+(n-1) \delta}\right)-q^{\frac{1}{2}} \Psi\left(F_{\alpha_1+(n-1) \delta}\right) y_2^{-1} x_2 \\
	&= & \frac{(-1)^n}{\left(q^{\frac{1}{2}}-q^{-\frac{1}{2}}\right)^2}\left(q^{\frac{-n^2-2 n+3}{2}} y_2^{-1} x_2 y_1^{-n} y_2 x_{-n}\right.  +q^{\frac{-n^2+6 n-1}{2}} y_2^{-1} x_2 y_1^{-n+4} y_2 x_{4-n} \\
	&& \left.-q^{\frac{-n^2-2 n+1}{2}} y_1^{-n} y_2 x_{-n} y_2^{-1} x_2 -q^{\frac{-n^2+6 n-3}{2}} y_1^{-n+4} y_2 x_{4-n} y_2^{-1} x_2\right) \\
	&= & \frac{(-1)^n}{\left(q^{\frac{1}{2}}-q^{-\frac{1}{2}}\right)^2}\left(q^{\frac{-n^2+2 n+1}{2}} y_1^{-n} x_2 x_{-n}+q^{\frac{-n^2+10 n-19}{2}} y_1^{-n+4} x_2 x_{4-n}\right. \\
	&& \left.-q^{-\frac{n^2+1}{2}} y_1^{-n} x_{-n} x_2-q^{\frac{-n^2+8 n-13}{2}} y_1^{-n+4} x_{4-n} x_2\right) \\
	&= & \frac{(-1)^{n+1}}{\left(q^{\frac{1}{2}}-q^{-\frac{1}{2}}\right)^2}\left(q^{-\frac{n^2+1}{2}} y_1^{-n}\left(x_{-n} x_2-q^{n+1} x_2 x_{-n}\right)\right. \\
	&& \left.+q^{\frac{3-(n-4)^2}{2}} y_1^{-n+4}\left(x_{4-n} x_2-q^{n-3} x_2 x_{4-n}\right)\right) .
\end{eqnarray*}
\endgroup
This proves (\ref{(F'nd)}) for any $n \ge 4$.

For $n=1$, note that $	x_{-1} x_2-q^2 x_2 x_{-1} =\left(1-q^2\right) X_\delta .$ By (\ref{(F'_delta)}), we obtain
\begin{eqnarray*}
\Psi \left( {{F_\delta }} \right) =  - \frac{{{q^{ - 1}}\left( {{q^2} - 1} \right)}}{{{{\left( {{q^{\frac{1}{2}}} - {q^{ - \frac{1}{2}}}} \right)}^2}}}y_1^{ - 1}{X_\delta } = \frac{{{q^{ - 1}}}}{{{{\left( {{q^{\frac{1}{2}}} - {q^{ - \frac{1}{2}}}} \right)}^2}}}y_1^{ - 1}\left( {{x_{ - 1}}{x_2} - {q^2}{x_2}{x_{ - 1}}} \right)
\end{eqnarray*}
For $n=2$, we have
\begin{eqnarray*}
	&&\Psi\left(F_{2 \delta}\right)=  q^{\frac{3}{2}} y_2^{-1} x_2 \Psi\left(F_{\alpha_1+\delta}\right)-q^{\frac{1}{2}} \Psi\left(F_{\alpha_1+\delta}\right) y_2^{-1} x_2 \\
	&= & \frac{1}{\left(q^{\frac{1}{2}}-q^{-\frac{1}{2}}\right)^2}\left(q^{-\frac{5}{2}} y_2^{-1} x_2 y_1^{-2} y_2 x_{-2}+q^2 y_2^{-1} x_2^2\right. \\
	&& \left.-q^{-\frac{7}{2}} y_1^{-2} y_2 x_{-2} y_2^{-1} x_2-q x_2 y_2^{-1} x_2\right) \\
	&= & \frac{1}{\left(q^{\frac{1}{2}}-q^{-\frac{1}{2}}\right)^2}\left(q^{\frac{1}{2}} y_1^{-2} x_2 x_{-2}-q^{-\frac{5}{2}} y_1^{-2} x_{-2} x_2+q^2 y_2^{-1} x_2^2-q^2 y_2^{-1} x_2^2\right) \\
	&= & -\frac{q^{-\frac{5}{2}}}{\left(q^{\frac{1}{2}}-q^{-\frac{1}{2}}\right)^2} y_1^{-2}\left(x_{-2} x_2-q^3 x_2 x_{-2}\right) .
\end{eqnarray*}
For $n=3$, we obtain
\begin{eqnarray*}
	&&\Psi\left(F_{3 \delta}\right)=  q^{\frac{3}{2}} y_2^{-1} x_2 \Psi\left(F_{\alpha_1+2 \delta}\right)-q^{\frac{1}{2}} \Psi\left(F_{\alpha_1+2 \delta}\right) y_2^{-1} x_2 \\
	&= & -\frac{1}{\left(q^{\frac{1}{2}}-q^{-\frac{1}{2}}\right)^2}\left(q^{-6} y_2^{-1} x_2 y_1^{-3} y_2 x_{-3}+q^{\frac{5}{2}} y_2^{-1} x_2 y_2 x_1\right. \\
	&& \left.-q^{-7} y_1^{-3} y_2 x_{-3} y_2^{-1} x_2-q^{\frac{3}{2}} y_2 x_1 y_2^{-1} x_2\right) \\
	&= & -\frac{1}{\left(q^{\frac{1}{2}}-q^{-\frac{1}{2}}\right)^2}\left(q^{-1} y_1^{-3} x_2 x_{-3}-q^{-5} y_1^{-3} x_{-3} x_2+q^{\frac{3}{2}}\left(x_2 x_1-x_1 x_2\right)\right) \\
	&= & \frac{q^{-5}}{\left(q^{\frac{1}{2}}-q^{-\frac{1}{2}}\right)^2} y_1^{-3}\left(x_{-3} x_2-q^4 x_2 x_{-3}\right) .
\end{eqnarray*}
The proof is complete.
\end{proof}

Similarly, we can rewrite $\Psi \left( {{F_{n\delta }}} \right)$ in terms of $S_n$ $(n\geq 0)$.
\begin{corollary}\label{5.7}
	For any $n\ge 1$, we have
	\begin{eqnarray*}
		\Psi \left( {{F_{n\delta }}} \right)	= \frac{{{{( - 1)}^n}{q^{ - \frac{{{n^2} + 1}}{2}}}}}{{{{\left( {{q^{\frac{1}{2}}} - {q^{ - \frac{1}{2}}}} \right)}^2}}}y_1^{ - n}\sum\limits_{j = 0}^{\left\lfloor {\frac{n}{2}} \right\rfloor }{\left( {{q^{n - j + 1}} - {q^j}} \right)\left( {1 + {\varepsilon _j}{q^{4n - 4}}y_1^2y_2^{ - 2}} \right){X^{(0,0,j,j)}}{S_{n - 2j}}}.
	\end{eqnarray*}
\end{corollary}

\begin{proof}
	By Theorem \ref{3.1} (1) and the bar-involution,
	for any $m,n\ge 0$, we have
	\begin{eqnarray}\label{(5.13)}
		&&x_{-(m+n+2)}x_{-m}-x_{-m}x_{-(m+n+2)}\nonumber\\
		&=&\left(q^{-\frac12}-q^{n+\frac12}\right)S_nX^{(0,0,m,m+1)}+x_{-(m+n+1)}x_{-(m+1)}-x_{-(m+1)}x_{-(m+n+1)}.
	\end{eqnarray}
	
	By repeatedly applying (\ref{(5.13)}),
	for any $n\ge 4$, we obtain
	\begin{eqnarray*}
		&&x_{-(n-2)}x_0-x_0x_{-(n-2)}\\
		&=&\left(q^{-\frac12}-q^{n-\frac72}\right)S_{n-4}X^{(0,0,0,1)}+x_{-(n-3)}x_{-1}-x_{-1}x_{-(n-3)}\\
		&\mathop =\limits^{(\ref{(5.13)})}&\sum\limits_{j=0}^{\left\lfloor\frac{n-4}{2}\right\rfloor}\left(q^{-\frac12}-q^{n-2j-\frac72}\right)S_{n-4-2j}X^{(0,0,j,j+1)}\\
		&=&\sum\limits_{j=2}^{\left\lfloor\frac n2\right\rfloor}\left(q^{-\frac12}-q^{n-2j+\frac12}\right)S_{n-2j}X^{(0,0,j-2,j-1)}.
	\end{eqnarray*}
	
	By Theorem \ref{3.1} (4), (5) and
	the bar-involution, for any $n\ge 2$, we obtain
	\begin{eqnarray*}
	&&x_{-n}x_2-q^{n+1}x_2x_{-n}\\
	&=&\left(1-q^{n+1}\right)S_n+q\left(1-q^{n-1}\right)X^{(0,0,1,1)}S_{n-2}\\
		&&+q^{n+\frac32}y_1X^{(0,0,1,1)}\left(x_{-(n-2)}x_0-x_0x_{-(n-2)}\right)\\
		&=&\left(1-q^{n+1}\right)S_n+q\left(1-q^{n-1}\right)X^{(0,0,1,1)}S_{n-2}\\
		&&+q^{n+\frac32}y_1X^{(0,0,1,1)}\sum\limits_{j=2}^{\left\lfloor\frac n2\right\rfloor}\left(q^{-\frac12}-q^{n-2j+\frac12}\right)S_{n-2j}X^{(0,0,j-2,j-1)}\\
		&=&\left(1-q^{n+1}\right)S_n+q\left(1-q^{n-1}\right)X^{(0,0,1,1)}S_{n-2}
		+\sum\limits_{j=2}^{\left\lfloor\frac n2\right\rfloor}\left(q^j-q^{n-j+1}\right)X^{(0,0,j,j)}S_{n-2j}\\
		&=&\sum\limits_{j=0}^{\left\lfloor\frac n2\right\rfloor}\left(q^j-q^{n-j+1}\right)X^{(0,0,j,j)}S_{n-2j}.
	\end{eqnarray*}	
\begingroup
	\interdisplaylinepenalty=0
	For $n=1$, we have $x_{-1}x_2-q^2x_2x_{-1}=(1-q^2)S_1.$ For $n=0$, by Theorem \ref{3.1} (3) and the bar-involution, we obtain $x_0x_2-qx_2x_0=1-q.$
	Therefore, for any $n\ge 0$, we have
	\begin{eqnarray}\label{(5.14)}
		x_{-n}x_2-q^{n+1}x_2x_{-n}=\sum\limits_{j=0}^{\left\lfloor\frac n2\right\rfloor}\left(q^j-q^{n-j+1}\right)X^{(0,0,j,j)}S_{n-2j}.
	\end{eqnarray}
	By Theorem \ref{5.6} (3) and (\ref{(5.14)}), for any $n\ge 4$, we obtain
\begin{eqnarray*}
		&&\Psi\left({{F_{n\delta}}}\right)=\frac{(-1)^{n+1}}{\left(q^{\frac12}-q^{-\frac12}\right)^2}\left(q^{-\frac{n^2+1}{2}}y_1^{-n}\sum\limits_{j=0}^{\left\lfloor\frac n2\right\rfloor}\left(q^j-q^{n-j+1}\right)X^{(0,0,j,j)}S_{n-2j}\right.\\
		&&\left.+q^{\frac{3-(n-4)^2}{2}}y_1^{-n+4}\sum\limits_{j=0}^{\left\lfloor\frac{n-4}{2}\right\rfloor}\left(q^j-q^{n-j-3}\right)X^{(0,0,j,j)}S_{n-4-2j}\right)\\
		&=&\frac{(-1)^nq^{-\frac{n^2+1}{2}}}{\left(q^{\frac12}-q^{-\frac12}\right)^2}y_1^{-n}\left(\sum\limits_{j=0}^{\left\lfloor\frac n2\right\rfloor}\left(q^{n-j+1}-q^j\right)X^{(0,0,j,j)}S_{n-2j}\right.\\
		&&\left.+q^{4n-6}y_1^4\sum\limits_{t=2}^{\left\lfloor\frac n2\right\rfloor}\left(q^{n-t-1}-q^{t-2}\right)X^{(0,0,t-2,t-2)}S_{n-2t}\right).\\
		&=&\frac{(-1)^nq^{-\frac{n^2+1}{2}}}{\left(q^{\frac12}-q^{-\frac12}\right)^2}y_1^{-n}\left(\sum\limits_{j=0}^{\left\lfloor\frac n2\right\rfloor}\left(q^{n-j+1}-q^j\right)X^{(0,0,j,j)}S_{n-2j}\right.\\
		&&\left.+q^{4n-6}\sum\limits_{t=2}^{\left\lfloor\frac n2\right\rfloor}q^4\left(q^{n-t-1}-q^{t-2}\right)y_1^2y_2^{-2}X^{(0,0,t,t)}S_{n-2t}\right)\\
		&=&\frac{(-1)^nq^{-\frac{n^2+1}{2}}}{\left(q^{\frac12}-q^{-\frac12}\right)^2}y_1^{-n}\left(\sum\limits_{j=0}^{\left\lfloor\frac n2\right\rfloor}\left(q^{n-j+1}-q^j\right)X^{(0,0,j,j)}S_{n-2j}\right.\\
		&&\left.+\sum\limits_{t=2}^{\left\lfloor\frac n2\right\rfloor}q^{4n-4}\left(q^{n-t+1}-q^t\right)y_1^2y_2^{-2}X^{(0,0,t,t)}S_{n-2t}\right)\\
		&=&\frac{(-1)^nq^{-\frac{n^2+1}{2}}}{\left(q^{\frac12}-q^{-\frac12}\right)^2}y_1^{-n}\sum\limits_{j=0}^{\left\lfloor\frac n2\right\rfloor}\left(q^{n-j+1}-q^j\right)\left(1+\varepsilon_jq^{4n-4}y_1^2y_2^{-2}\right)X^{(0,0,j,j)}S_{n-2j}.
	\end{eqnarray*}
	\endgroup
	
For $n=1,2,3$, note that $R'(n)=0$. Thus, by (\ref{(5.14)}), we have
	\[\Psi\left({{F_{n\delta}}}\right)=\frac{(-1)^nq^{-\frac{n^2+1}{2}}}{\left(q^{\frac12}-q^{-\frac12}\right)^2}y_1^{-n}\sum\limits_{j=0}^{\left\lfloor\frac n2\right\rfloor}\left(q^{n-j+1}-q^j\right)X^{(0,0,j,j)}S_{n-2j}.\]
This completes the proof.
\end{proof}

\section*{Acknowledgments}
This work was supported by NSF of China (No. 12371036).
\section*{Declaration of AI use}
 ChatGPT was used  for checking proofs and improving the clarity of the manuscript. All mathematical results, arguments and proofs were carefully developed and verified by the authors.


\begin{thebibliography}{99}

\bibitem{BR}
Berenstein A, Rupel D. Quantum cluster characters of Hall algebras. Selecta Math (NS), 2015, 21: 1121--1176

\bibitem{AA}
Berenstein A, Zelevinsky A. Quantum cluster algebras. Adv Math, 2005, 195: 405--455

\bibitem{CLW21a}
Chen X, Lu M, Wang W. A Serre presentation for the $\imath$quantum groups. Transform Groups, 2021, 26: 827--857

\bibitem{Da93}
Damiani I. A basis of type Poincar\'e--Birkhoff--Witt for the quantum algebra of $\widehat{\mathfrak{sl}}_2$. J Algebra, 1993, 161: 291--310

\bibitem{DXC23}
Ding M, Xu F, Chen X. Recursive formulas for the Kronecker quantum cluster algebra with principal coefficients. Sci China Math, 2023, 66: 1933--1948

\bibitem{Dri85}
Drinfeld V. Hopf algebras and quantum Yang-Baxter equation. Soviet Math Dokl, 1985, 32: 254--258

\bibitem{ca1}
Fomin S, Zelevinsky A. Cluster algebras. I. Foundations. J Amer Math Soc, 2002, 15: 497--529

\bibitem{ca2}
Fomin S, Zelevinsky A. Cluster algebras. II. Finite type classification. Invent Math, 2003, 154: 63--121

\bibitem{FZ25}
Fu C, Zhang H. From quantum groups to quantum cluster algebras. arXiv:2511.14240v1, 2025

\bibitem{Goff24}
Goff O. The $q$-Onsager algebra and the quantum torus. J Combin Theory Ser A, 2024, 208: Paper No. 105939, 22 pp

\bibitem{Goff25}
Goff O. Using the quantum torus to investigate the $q$-Onsager algebra. arXiv:2504.13362v1, 2025

\bibitem{HJY}
Huang J, Chen X, Ding M, Xu F. Fundamental relations in quantum cluster algebras. arXiv:2509.11635v1, 2025

\bibitem{Ip18}
Ip I. Cluster realization of $U_q(\mathfrak{g})$ and factorizations of the universal $R$-matrix. Selecta Math (NS), 2018, 24: 4461--4553

\bibitem{Jim85}
Jimbo M. A $q$-difference analogue of $U(g)$ and the Yang-Baxter equation. Lett Math Phys, 1985, 10: 63--69

\bibitem{Ko14}
Kolb S. Quantum symmetric Kac-Moody pairs. Adv Math, 2014, 267: 395--469

\bibitem{Let99}
Letzter G. Symmetric pairs for quantized enveloping algebras. J Algebra, 1999, 220: 729--767


\bibitem{LRZ25}
Lu M, Ruan S, Zhang H. Analogue of Feigin's map on $\imath$quantum groups of split type. Int Math Res Not, 2025, 2025(19): rnaf297

\bibitem{LW22}
Lu M, Wang W. Hall algebras and quantum symmetric pairs I: foundations. Proc Lond Math Soc, 2022, 124: 1--82

\bibitem{LW23}
Lu M, Wang W. Hall algebras and quantum symmetric pairs of Kac-Moody type. Adv Math, 2023, 430: Paper No. 109215, 56 pp

\bibitem{G1}
Lusztig G. Canonical bases arising from quantized enveloping algebras. J Amer Math Soc, 1990, 3: 447--498

\bibitem{G2}
Lusztig G. Quivers, perverse sheaves, and quantized enveloping algebras. J Amer Math Soc, 1991, 4: 365--421

\bibitem{G3}
Lusztig G. Total positivity in reductive groups. In: Lie Theory and Geometry. Progr Math, vol. 123. Boston: Birkh\"auser, 1994: 531--568

\bibitem{SS19}
Schrader G, Shapiro A. A cluster realization of $U_q(\mathfrak{sl}_n)$ from quantum character varieties. Invent Math, 2019, 216: 799--846

\bibitem{Shen}
Shen L. Cluster nature of quantum groups. arXiv:2209.06258, 2022

\bibitem{So24}
Song J. Cluster realisations of $\imath$quantum groups of type AI. Proc Lond Math Soc, 2024, 129: e12636

\bibitem{Ter19}
Terwilliger P. Using Catalan words and a $q$-shuffle algebra to describe a PBW basis for the positive part of $U_q(\widehat{\mathfrak{sl}}_2)$. J Algebra, 2019, 525: 359--373

\end{thebibliography}
\end{document}